\documentclass[smallextended]{svjour3}
\smartqed  

\makeatletter
\def\cl@chapter{\@elt {theorem}}
\makeatother
\usepackage{graphicx}
\usepackage{amsmath,amssymb}
\usepackage{graphicx}
\usepackage{mhchem}
\usepackage[capitalise]{cleveref}
\usepackage{xcolor}
\usepackage{verbatim}
\usepackage[misc,geometry]{ifsym} 

\usepackage{algorithm}
\usepackage{algorithmic}
\usepackage{url}

\usepackage[font=small,skip=0pt]{caption}  % to reduce space between caption and figure,table
\newcommand{\COMMA}{ \; ,}
\newcommand{\PERIOD}{ \; .}

\newcommand{\REVa}[1]{{\color{black}#1}}
\newcommand{\REVb}[1]{{\color{black}#1}}

\def\CC{{C\nolinebreak[4]\hspace{-.05em}\raisebox{.4ex}{\tiny\bf ++}}}

\begin{document}

\title{ $S$-Leaping: An adaptive, accelerated stochastic simulation algorithm, bridging $\tau$-leaping and $R$-leaping}

\titlerunning{$S$-Leaping: bridging $\tau$-leaping and $R$-leaping}        % if too long for running head

\author{
Jana Lipkov\'{a} \and
Georgios Arampatzis  \and
Philippe Chatelain  \and
Bjoern Menze \and
Petros Koumoutsakos
}


\institute{
J. Lipkov\'{a} and B. Menze \at Department of Informatics, Technical University of Munich, DE-85748, Germany
\and
P. Chatelain \at Institute of Mechanics, Materials and Civil Engineering, Universit\'{e} catholique de Louvain, 1348 Louvain-la-Neuve, Belgium
\and
G. Arampatzis and P. Koumoutsakos (\Letter  \hspace{0pt} petros@ethz.ch)  \at Computational Science and Engineering Laboratory, ETH Zurich, Zurich, CH-8092, Switzerland
}

\date{Received: date / Accepted: date}

\maketitle

%==================================================
%==========	 Abstract
%==================================================
\begin{abstract}
We propose the $S$-leaping algorithm for the acceleration of Gillespie's stochastic simulation algorithm that combines the advantages of the two main accelerated methods; the $\tau$-leaping and \mbox{$R$-leaping} algorithms. These algorithms are known to be efficient under different conditions; the  $\tau$-leaping is efficient for non-stiff systems or systems with partial equilibrium, while the $R$-leaping performs better in stiff system thanks to an efficient sampling procedure. However, even a small change in a system's set up can critically affect the nature of the simulated system and thus reduce the efficiency of an accelerated algorithm. The proposed algorithm combines the efficient time step selection from the $\tau$-leaping with the effective sampling procedure from the $R$-leaping algorithm. The $S$-leaping is shown to  maintain its efficiency under different conditions and in the case of large and stiff systems or systems with fast dynamics, the $S$-leaping outperforms both methods. We demonstrate the performance and the accuracy of the $S$-leaping in comparison with the $\tau$-leaping and $R$-leaping on a number of benchmark systems involving biological reaction networks.


\end{abstract}
%============================================

%==================================================
%==========	Section: Introduction
%==================================================
\section{Introduction}
\label{sec:intro}

The celebrated Gillespie's stochastic simulation algorithm (SSA) \cite{Gillespie1976,Gillespie:1977} simulates continuous-time Markov chains systems. An example of such system is a well-stirred chemically reacting system with small population of reactants \cite{Anderson:2011}.
The SSA is an exact numerical algorithm. However, since SSA allows to simulate only one reaction event per time step, it becomes computationally costly for large systems and long time scales. Over the years, several algorithms were proposed to accelerate the SSA at the expense of sacrificing its accuracy. The most prominent are the $\tau$-leaping \cite{Gillespie:2001} with its  further enhancements \cite{Cao1:2005,Cao:2006,Cao:2005,Cao:2007,Rathinam:2003,Tian:2004} and the  $R$-leaping algorithm \cite{Auger:2006,Mjolsness:2009}. Other accelerated algorithms involve the FLAVOR-SSA, where flow averaging is used to accelerate the simulation \cite{Bayati:2010d},  coupling of multi-scale frameworks with any stochastic simulation algorithm \cite{Koumoutsakos:2013a} and an adaptive mesh refinement algorithm for reaction-diffusion systems \cite{Bayati:2011}.
One can finally mention a special class of algorithms which achieve both exact, SSA-like, sampling of the reaction events and computational acceleration, as initiated by the Exact $R$-leaping~\cite{Mjolsness:2009}; the acceleration offered by such techniques is however weaker than in the $\tau$-leaping and $R$-leaping algorithms.

The $\tau$-leaping algorithm \cite{Gillespie:2001} accelerates the SSA by advancing the state of the system by a larger time step $\tau$, allowing multiple reaction events to occur within the preselected time step. The number of firings of each reaction channel at each time step is a random variable that follows Poisson distribution. On the other hand, the $R$-leaping algorithm preselects the total number of reaction firings $L$ \cite{Auger:2006}. The time step needed for those $L$ reactions events to occur follows a Gamma distribution and the number of firings of each reaction follows a multinomial distribution, which can be efficiently sampled through correlated binomial distributions. Both approximate algorithms are valid under the \textit{leap condition} which states that the propensities must remain approximately constant during each simulation step.

Each of these algorithms is efficient under different conditions. In non-stiff systems, the  $\tau$-leaping is more effective than the $R$-leaping algorithm. In addition, the implicit extension of the $\tau$-leaping for stiff systems where some reaction channels appear in partial equilibrium \cite{Cao:2007}, allows to advance the system with bigger time steps, which yields to significant speed-up over the explicit $R$-leaping method. However, the sampling procedure in the $\tau$-leaping method requires to draw one random number for each reaction channel.  This is especially inefficient in big and stiff systems, where only few reaction channels are fired per time step.
On the other hand, since the samples in the $R$-leaping are drawn from a correlated probability distribution, the amount of drawn random numbers can be reduced by reordering the reaction indices in a way that the most probable reaction channels are sampled first. This yields appreciable computational savings in big and stiff systems. %The weakness of $R$-leaping is that it does not incorporate the speed-up provided by a partial equilibrium of some reaction channels, since the number of firings is preselected and the time step is computed afterwards.

In this paper we present the $S$-leaping algorithm as an efficient coupling of both methods. Our algorithm uses the efficient time step selection procedure present in the $\tau$-leaping. This feature allows the $S$-leaping to exploit the advantage of implicit formulation for stiff systems with partial equilibrium. In addition, the $S$-leaping estimates the total number of firings within a preselected time interval as a sample from Poisson distribution. The knowledge of the total number of reaction channels allows to draw individual firings from the correlated binomial distributions, with further optimization through reordering of reaction channels in big and stiff systems. Thus, the $S$-leaping algorithm provides an effective fusion of both methods. The name of the method was chosen so that it represents the position of the $S$-leaping between the $R$-leaping and $\tau$-leaping method.

The paper is organized as follows. In Section \ref{sec:background} we provide a brief specification of the SSA, $\tau$-leaping and $R$-leaping algorithms.
The $S$-leaping algorithm is introduced in Section \ref{sec:S}.
In Section \ref{sec:numerics1}, the $S$-leaping method is tested on four benchmark cases, a non-stiff, a stiff, a fast dynamics and a large reaction network.  We conclude with a summary in Section \ref{sec:summary}. 
%==================================================

%==================================================
%==========		Section: Background
%==================================================
\section{Background} \label{sec:background}

We consider a well-stirred system that contains $N$ molecular species $\{ S_{1},\dots,$ $S_{N} \}$ that can react through $M$ chemical reactions channels $\left\{ R_{1}, \dots, R_{M}  \right\}$. In what follows, the letter $i \in \{1, \dots, N\}$ will be used for chemical species, e.g., $S_{i}$, and the letter $j \in \{1, \dots, M\}$ for chemical reactions, e.g., $R_{j}$.
The state of the system is characterized by the \textit{state vector} $\vec{X}(t) = \left ( X_{1}(t),  \ldots,  X_{N}(t)\right )$, where $X_{i}(t)$ denotes the number of molecules $S_{i}$ at time $t$. The dynamics of each reaction channel $R_{j}$ are being  characterized by a \textit{propensity function} $a_{j}$ and a \textit{state change vector} $\vec{\nu}_{j}=({\nu}_{1j},\ldots,{\nu}_{Nj})$. Given the state vector $\vec{x}=(x_1,\ldots,x_N)=\textbf{X}(t)$, the quantity $a_{j}(\vec{x})dt$ gives the probability that the reaction $R_{j}$ will occur  in the next infinitesimal time interval $\left[ t, t+dt \right)$. The state change vector $\vec{\nu}_j$ gives the change in the molecular population caused by one reaction $R_{j}$. \REVa{Finally, we define $a_{0}$ as the sum of all propensity functions $a_0(\vec{x})=\sum_{j=1}^{M} a_j(\vec{x})$.} 

For the rest of the paper we will use the notation $\mathcal{B}, \Gamma, \mathcal{E}, \mathcal{M}$, and $\mathcal{P}$ to denote the probability distribution function of the binomial, the gamma, the exponential, the multinomial and the Poisson distribution, respectively. The same notation will be used to denote the function that produces pseudo-random numbers from the respective distribution.  With $\lfloor x \rceil$ we will denote the closest integer to $x$. %\REVa{Finally,  $a_0$ is the function defined through  $a_0(\vec{x})=\sum_{j=1}^{M} a_j(\vec{x}) $}.

%==============================
%       	 Subsection: SSA
%==============================
%
%
\subsection{The Stochastic Simulation Algorithm}
The SSA \cite{Gillespie:1977}  is an exact algorithm for simulating the time evolution of well-stirred chemically reacting systems.  It is an exact algorithm in the sense that the generated sample paths are distributed according to the solution of the corresponding chemical master equation \cite{Gillespie:1977}. However, since SSA simulates only one reaction event per time step, it becomes inefficient for most realistic systems. The SSA algorithm is summarised in \cref{alg:SSA}.
%~~~~~~~~~~~~~~~~~~~~~~~~~~~
\begin{algorithm}
\caption{\label{alg:SSA} Stochastic Simulation Algorithm (SSA)}
\algsetup{indent=1.5em}
{
\fontsize{8}{12}\selectfont
\begin{algorithmic}[1]
\STATE Initialise: $T_{\textrm{end}}$, $\vec{x} \gets \vec{X}(0)$, $t \gets 0$
\WHILE{$ t<T_{\textrm{end}}$}
\STATE Compute $a_{j}(\vec{x})$ for $ j= 1, \dots , M$ and $a_{0}(\vec{x})$%\gets \sum_{j=1}^{M} a_{j}( \vec{x} )$
\STATE $\tau \gets \mathcal{E}\left( 1/a_{0}(\vec{x})\right)$
\STATE Choose the $j$-th reaction with probability ${a_{j}(\vec{x})} / { a_{0}(\vec{x} )} $
\STATE $\vec{x} \gets \vec{x} + \vec{\nu}_{j}$
\STATE $t \gets t+\tau$
\ENDWHILE
\end{algorithmic}
}
\end{algorithm}
%~~~~~~~~~~~~~~~~~~~~~~~~~~~
%
%
%==============================

%==============================
%       	 Subsection: Leap Method
%==============================
\subsection{Approximate accelerated stochastic simulation algorithms}
Several approximate stochastic simulation algorithms \cite{Auger:2006,Cao:2006,Cao:2005,Gillespie:2001} have been introduced to accelerate the SSA by advancing the system with larger time steps, allowing to fire more reactions per time step. The accurate advancement of the system is limited by the so called \textit{leap condition}, which states that propensities $a_{j}(\vec{x})$ should remain approximately constant over the time interval $[ t , t+\tau)$,
\begin{equation}\label{eq:LeapCondition}
| a_{j}(\textbf{X}(t+\tau)) - a_{j}(\textbf{X}(t)) | \leq \varepsilon \, a_{0}(\vec{x}), \quad j=  1,\ldots,M \COMMA
\end{equation}
where $0 < \varepsilon \ll1$ is a user defined parameter that controls the models accuracy.
%==============================

%==============================
%  SubSub: NN-TAU
%==============================
\subsubsection{Non-negative $\tau$-leaping}
The $\tau$-leaping algorithm \cite{Gillespie:2001} pre-selects a deterministic time step $\tau$, much bigger that the mean stochastic time step of SSA. Then, the number of times $k_{j}^{\mathcal{P}}$ the reaction $R_j$ will be fired during the time interval $[t,t+\tau)$ is sampled from a Poisson distribution with parameter $a_{j}(\vec{x})\tau$. Since the Poisson random variables $k_{j}^{\mathcal{P}}$ are unbounded, the algorithm might result in negative populations. To overcome this problem a non-negative version of the $\tau$-leaping algorithm was proposed in \cite{Cao:2005}. The algorithm identifies the critical reactions, those which are $N_c$ firings from exhausting one of its reactants. No more than one critical reaction can occur within the time leap $\tau$, while multiple non-critical reactions are allowed. The critical reaction is handled by the SSA, while the non-critical reactions are modelled by the $\tau$-leaping method. Several methods \cite{Cao:2006,Gillespie:2001,Gillespie:2003} were introduced for the computation of the leap length $\tau$. The most efficient one \cite{Cao:2006} selects $\tau$ by
%--------------------------------------------------------
\begin{equation}\label{eq:ExpTauTimeStep}
\tau = \min_{i\in I_{\textrm{rs}}} \left\{    \frac{ \max \left\{  \frac{\varepsilon x_{i}}{g_{i}(\vec{x})} , 1\right\}  }{ \lvert \, \mu_{i} (\vec{x}) \, \rvert },   \frac{  \max\left\{   \frac{\varepsilon x_{i}}{ g_{i}(\vec{x})} , 1\right\}^2   }{  \lvert \, \sigma_{i}^{2}(\vec{x}) \,   \rvert  } \right\} \COMMA
\end{equation}
%--------------------------------------------------------
for $\vec{x}=\vec{X}(t)$ and $I_{\textrm{rs}}$  the set of indices of all reactant species. The factor $g_{i}$ takes into account the highest order of reaction, denoted as $h_i$, in which species $S_{i}$ appears as a reactant,
%--------------------------------------------------------
\begin{equation}\label{eq:g_i}
g_{i}(\vec{x})= h_i  + \frac{h_i}{n_i} \sum_{j=1}^{n_i - 1} \frac{j}{x_i - j} \COMMA
\end{equation}
%--------------------------------------------------------
where $n_i$ denotes the maximum number of $S_{i}$ molecules required by any of the highest order reactions \cite{Sandmann:2009}. Finally, the terms $\mu_{i}$ and $\sigma^{2}_{i}$ are given by
%--------------------------------------------------------
\begin{eqnarray}
\mu_{i} (\vec{x} ) &=&  \sum_{j \in J_{\textrm{ncr}}} \nu_{ij} \, a_{j}(\vec{x}), \quad \forall i \in I_{\textrm{rs}} \label{eq:mu} \COMMA \\
\sigma_{i}^{2}(\vec{x} ) &=& \sum_{j \in J_{\textrm{ncr}}} \nu_{ij}^{2} \, a_{j} (\vec{x}), \quad \forall i \in I_{\textrm{rs}}, \label{eq:sigma} \COMMA
\end{eqnarray}
%--------------------------------------------------------
where $J_{\textrm{ncr}}$ is the set of all non-critical reactions. The non-negative $\tau$-leaping algorithm is outlined in \cref{alg:Tau}.
%
%
%~~~~~~~~~~~~~~~~~~~~~~~~~~~
\begin{algorithm}
\algsetup{indent=1.5em}
\caption{\label{alg:Tau}Non-negative $\tau$-leaping}
{
\fontsize{8}{11}\selectfont
\begin{algorithmic}[1]
\STATE Initialise: $T_{\textrm{end}}$, $\vec{x} \gets \vec{X}(0)$, $t \gets 0$, $N_{\textrm{c}}\gets10$.
\WHILE{$t<T_{\textrm{end}}$}
\STATE Compute $a_{j}(\vec{x})$ for $j=1, \dots , M$ and $a_{0}(\vec{x})$% \gets \sum_{j=1}^{M} a_{j}( \vec{x})$.
\STATE Compute the list of critical reactions $J_\textrm{crit}$. The reaction $R_{j}$ is critical if:
$$
a_{j}(\vec{x})>0 \quad \textrm{and} \quad \min_{i}  \left \lfloor \frac{ x_{i} }{ |\nu_{ij}| } \right \rceil \leq N_c
$$
\STATE Compute time the step $\tau_{1}$ by \cref{eq:ExpTauTimeStep}
\IF{$\tau_{1} < 10 \frac{1}{  a_{0}(\vec{x})}$}\label{alg:tau:line:tau}
\STATE Execute 100 steps of the SSA
\ELSE
\STATE $a_{0}^{\textrm{c}}(\vec{x}) = \sum_{j \in J_\textrm{crit}} a_{j}(\vec{x})\quad$	and $\quad\tau_{2} \gets \mathcal{E}(1/a^{\textrm{c}}_{0}(\vec{x} ) )$  \COMMENT{time of critical reaction}
\IF{$\tau_{1} \leq \tau_{2}$}
\STATE $\tau \gets \tau_1$
\STATE $k_{j} \gets \mathcal{P} (a_{j}(\vec{x}) \tau), \quad j \notin J_\textrm{crit} $
\STATE $k_{j} = 0,\quad j \in J_\textrm{crit}$
\ELSE \STATE $\tau \gets \tau_2$
\STATE Choose $j_{\textrm{c}}$ with probability $ a_{j_{\textrm{c}}}(\vec{x}) / a_0^c(\vec{x})$ and $j_c\in J_\textrm{crit}$
\STATE $k_{j_c} \gets 1$
\STATE $k_{j} \gets 0$ for $ j\in J_\textrm{crit}$ and $j\neq j_{\textrm{c}}$
\STATE $k_{j} \gets \mathcal{P}(a_{j}(\vec{x})\tau)$ for $j \notin J_\textrm{crit}$
\ENDIF
\ENDIF
\IF{ there is a negative component in $\vec{x} + \sum_{j=1}^{M} k_{j}\vec{\nu}_{j}$ }
\STATE $\tau_{1} \gets \tau/2$ and go to \ref{alg:tau:line:tau}.
\ELSE{}
\STATE $\vec{x} \gets \vec{x} + \sum_{j=1}^{M} k_{j}\vec{\nu}_{j}$
\STATE $t \gets t+\tau$.
\ENDIF
\ENDWHILE
\end{algorithmic}
}
\end{algorithm}
%~~~~~~~~~~~~~~~~~~~~~~~~~~~
%
%==============================
%
%
%
%
%
%==============================
% SubSec: Adaptive Tau
%==============================

\subsubsection{Adaptive $\tau$-leaping}\label{subsec:AdaptiveTau}
An adaptive version of the $\tau$-leaping algorithm was introduced in \cite{Cao:2007}. It automatically alternates between the explicit (\cref{alg:Tau}) and implicit $\tau$-leaping \cite{Rathinam:2003} algorithm. The implicit $\tau$-leaping algorithm is inspired by the implicit Euler method for differential equations. Ideally, we would like to compute the state $\vec{X} (t+\tau)$ as
%-----------------------------------
\begin{equation}\label{eq:ImplTauUpdate}
\vec{X}(t + \tau) = \vec{X}(t) + \sum_{j=1}^{M} \vec{\nu}_{j} \, k_{j}^{\mathcal{P}}(\vec{X}(t+\tau)) \PERIOD
\end{equation}
%-----------------------------------
However, this would require the generation of random samples from a Poisson distribution with unknown parameter,
%-----------------------------------
\begin{eqnarray*}
k_{j}^{\mathcal{P}}(\vec{X}(t+\tau)) \sim \mathcal{P}(a_{j}\left( \vec{X}(t+\tau \right) \tau) \PERIOD
\end{eqnarray*}
%-----------------------------------
To avoid this difficulty, a partial implicit approach  was introduced in \cite{Rathinam:2003}. If $k^{\mathcal{P}}_{j}$ is a random variable that follows a  Poisson distribution with mean $a_{j}\tau$, then $k^{\mathcal{P}}_{j}$ can be expressed as a sum of a random variable with mean $a_{j}\tau$ and zero mean random variable $k^{\mathcal{P}}_{j} - a_{j}\tau$,
%-----------------------------------
\begin{equation}\label{eq:PoissSplit}
k^{\mathcal{P}}_{j} = a_{j} \tau+ k^{\mathcal{P}}_{j} - a_{j} \tau \PERIOD
\end{equation}
%-----------------------------------
The partial implicit approach evaluates the variable $a_{j}\tau$ at the state $\vec{X}(t+\tau)$ and the zero mean variable $k^{\mathcal{P}}_{j} - a_{j}\tau$ at the state $\vec{X}(t)$.  Applying this approach to the firings $k_{j}^{\mathcal{P}}$ in \cref{eq:ImplTauUpdate} leads to the following implicit system of equations,
%-----------------------------------
\begin{equation}\label{eq:tauIMPLSystem}
\vec{x}^\prime = \vec{x} + \sum_{j=1}^{M}\vec{\nu}_{j} a_{j}(\vec{x}^\prime) \tau  + \sum_{j=1}^{M}\vec{\nu}_{j}\left(  k_{j}^{\mathcal{P}}(\vec{x} )  - a_{j}(\vec{x}) \tau \right) \PERIOD
\end{equation}
%-----------------------------------
for $\vec{x}=\vec{X}(t)$ and $\vec{x}^\prime=\vec{X}(t+\tau)$.
If we denote by $\vec{X}^{\star}$ the solution of the above implicit system, which can be obtained with Newton-Raphson method, the implicit state update in \cref{eq:ImplTauUpdate} is given by,
%-----------------------------------
\begin{equation}
\vec{X}(t+\tau) = \vec{X}(t) + \sum_{j=1}^{M} \vec{\nu}_{j} k_{j}^{\mathcal{P}\star} ,
\end{equation}
%-----------------------------------
where
%-----------------------------------
\begin{equation}\label{eq:Tau:Imp:kj}
k_{j}^{\mathcal{P}\star} = \Big \lfloor  a_{j}(\vec{X}^{\star}) \tau + k_{j}^{\mathcal{P}}(\vec{X}(t)) - a_{j}(\vec{X}(t))\tau \Big\rceil \PERIOD
\end{equation}
%-----------------------------------
The rounding in \cref{eq:Tau:Imp:kj} ensures that the updated population will remain integer.

Implicit numerical methods provide an efficient way for solving stiff systems since they advance the system with  bigger time steps than explicit methods. While  implicit methods for differential equations  are unconditionally stable, the time step in the implicit leaping methods is bounded by the leap condition of \cref{eq:LeapCondition}.  The computation of the implicit leap step $\tau$ under the condition of partial equilibrium was introduced in \cite{Cao:2007}. The assumption is that if some reaction channels are in equilibrium or close to a partial equilibrium, then the net change of their propensities would be small. Thus the dynamics of the system would be driven by the reactions outside the equilibrium and the implicit time step can be computed as
%-----------------------------------
\begin{equation}\label{eq:ImplTau}
\tau^{(\textrm{im})} =\min_{i\in I_{\textrm{rs}}}  \left\{\frac{  \max \left\{  \frac{\varepsilon x_{i}}{g_{i}(\vec{x})} , 1\right\} }{ \lvert \mu_{i}^{(\textrm{im})} (\vec{x}) \rvert   },\frac{\max\left\{   \frac{\varepsilon x_{i}}{g_{i}(\vec{x})} ,1\right\}^2    }{   \sigma_{i}^{(\textrm{im})}(\vec{x})^2 }  \right\} \COMMA
\end{equation}
%-----------------------------------
where $g_i$ is given by \cref{eq:g_i} and  $\mu_{i}^{(\textrm{im})}$,  $\sigma_{i}^{(\textrm{im})}$ are given by
%-----------------------------------
\begin{eqnarray}
\mu_{i}^{(\textrm{im})}(\vec{x}) 		&=& \sum_{j \in J_{\textrm{necr}}} {\nu}_{ij} \, a_{j}(\vec{x}), \quad \forall i \in I_{\textrm{rs}} \COMMA \\
\sigma_{i}^{(\textrm{im})}(\vec{x})^2 	&=& \sum_{j \in J_{\textrm{necr}}} {\nu}_{ij}^{2} \, a_{j}(\vec{x}), \quad \forall i \in I_{\textrm{rs}} \COMMA
\end{eqnarray}
%-----------------------------------
for $\vec{x}=\vec{X}(t)$. Here,  $J_{\textrm{necr}}$ denotes the set of indices of the reaction channels that are neither critical nor in partial equilibrium.

In general, it is difficult to detect which reaction channels are currently in  partial equilibrium, however, it can be easily detected for reversible reactions \cite{Cao:2007}. Let $R_{+}$ and $R_{-}$ denote a pair of reversible reactions, with the corresponding propensity functions $a_{+}$ and $a_{-}$. If the reaction $R_{+}$ and $R_{-}$ are in partial equilibrium, their propensities must be similar,
%-----------------------------------
\begin{equation}\label{eq:PEC}
\lvert a_{+}(\vec{x}) - a_{-}(\vec{x}) \rvert \leq \delta \min \{ a_{+}(\vec{x}), a_{-}(\vec{x}) \} \COMMA
\end{equation}
%-----------------------------------
where $\delta$ is a small positive number, usually chosen  around 0.05 \cite{Cao:2007}. The adaptive $\tau$-leaping algorithm is outlined in \cref{alg:AdaptTau}.
%
%
%~~~~~~~~~~~~~~~~~~~~~~~~~~~
\begin{algorithm}
\caption{\label{alg:AdaptTau} Adaptive $\tau$-leaping}
\algsetup{indent=1.5em}
{
\fontsize{8}{12}\selectfont
\begin{algorithmic}[1]
\STATE Initialise: $T_{\textrm{end}}$, $\vec{x} \gets \textbf{X}(0)$, $t \gets 0$, $N_{\textrm{c}}\gets10$.
\WHILE{$t<T_{\textrm{end}}$}
	\STATE Compute $a_{j}(\vec{x})$ for $j=1, \dots , M$ and $a_{0}(\vec{x})$% \gets \sum_{j=1}^{M} a_{j}( \vec{x})$
        \STATE Compute the list of critical reactions $J_\textrm{crit}$. The reaction $R_{j}$ is critical if:
        $$
	a_{j}(\vec{x})>0 \quad \textrm{and} \quad \min_{i}  \left \lfloor \frac{ x_{i} }{ |\nu_{ij}| } \right \rceil \leq N_c
        $$
	\STATE Compute $\tau^{(\textrm{ex})}$ using \cref{eq:ExpTauTimeStep} and $\tau^{(\textrm{im})}$ using \cref{eq:ImplTau}
	\IF{$ \tau^{(\textrm{im})} > 100\; \tau^{(\textrm{ex})} $}
		\STATE The system is stiff and $\tau_{1} \gets \tau^{(\textrm{im})}$
	\ELSE{}
		\STATE The system is non-stiff and $\tau_{1}\gets \tau^{(\textrm{ex})}$
	\ENDIF
	\IF{$ \tau_{1} \leq 10 \frac{1}{a_{0}(\vec{x})} $} \label{alg:adatau:line:tau}
		\STATE Execute 100 steps of the SSA.
	\ELSE
	\STATE $a_{0}^{\textrm{c}}(\vec{x}) = \sum_{j \in J_\textrm{crit}} a_{j}(\vec{x})\quad$	and $\quad\tau_{2} \gets \mathcal{E}(1/a^{\textrm{c}}_{0}(\vec{x} ) )$  \COMMENT{time of critical reaction}
		\IF{$ \tau_2 > \tau_1$}
			\STATE $\tau \gets \tau_1$
			\IF{ the system is currently stiff}
				\STATE Compute $k_{j}$ using \cref{eq:Tau:Imp:kj} for $j\notin J_\textrm{crit}$
			\ELSE
        			\STATE $k_{j} \gets \mathcal{P}(a_{j}(\vec{x}))$ for $j\notin J_\textrm{crit}$
			\ENDIF
		\ELSE
			\STATE $\tau \gets \tau_2$
			\STATE Choose $j_{\textrm{c}}$ with probability $ a_{j_{\textrm{c}}}(\vec{x}) / a_0^c(\vec{x})$ and $j_c\in J_\textrm{crit}$
			\STATE $k_{j_c} \gets 1$
			\STATE $k_{j} \gets 0$ for $ j\in J_\textrm{crit}$ and $j\neq j_{\textrm{c}}$
			\IF{ $ \tau_2 < \tau^{ex}$ \OR the system is non-stiff }
				\STATE $k_{j} \gets \mathcal{P}(a_{j}(\vec{x})\tau)$ for $j\notin J_\textrm{crit}$
				\ELSE
				\STATE Compute $k_{j}$ using \cref{eq:Tau:Imp:kj} for $j\notin J_\textrm{crit}$
			\ENDIF
	\ENDIF
	\IF{ there is a negative component in $\vec{x} + \sum_{j=1}^{M} k_{j}\vec{\nu}_{j}$ }
		\STATE $\tau_{1} \gets \tau/2$ and go to \ref{alg:adatau:line:tau}
		\ELSE{}
		\STATE  $\vec{x} \gets \vec{x} + \sum_{j=1}^{M} k_{j}\vec{\nu}_{j}$
		\STATE $t \gets t+\tau$
	\ENDIF
\ENDIF
\ENDWHILE
\end{algorithmic}
}
\end{algorithm}
%%~~~~~~~~~~~~~~~~~~~~~~~~~~~
%
%
%
%==============================
%
%
%
%
%
%
%==============================
%  SubSub: R
%==============================
%
\subsubsection{R-leaping}\label{subsec:R}
The $R$-leaping algorithm \cite{Auger:2006},  instead of prescribing the time-step, it imposes the total number of reactions $L$ that can be fired during the next time interval. Under the leap condition of \cref{{eq:LeapCondition}},  the number of firings is computed as \cite{Auger:2006},
%--------------------------------------------------------
\begin{equation}\label{eq:R:L}
L = a_{0}(\vec{x})\min_{i\in I_{\textrm{rs}}} \left\{    \frac{ \max \left\{  \frac{\varepsilon x_{i}}{g_{i}(\vec{x})} , 1\right\}  }{ \lvert \, \mu_{i} (\vec{x}) \, \rvert },   \frac{  \max\left\{   \frac{\varepsilon x_{i}}{ g_{i}(\vec{x})} , 1\right\}^2   }{  \lvert \, \sigma_{i}^{2}(\vec{x}) \,  \rvert -  \rvert \mu_{i}^{2}(\vec{x})/a_{0}(\vec{x}) \rvert} \right\} \COMMA
\end{equation}
%--------------------------------------------------------
for $\vec{x}=\vec{X}(t)$, $I_{\textrm{rs}}$  the set of indices of all reactant species and the terms $g_{i}$, $\mu_{i}(\vec{x})$ and $\sigma^{2}_{i}(\vec{x})$ given by \cref{eq:g_i}, \eqref{eq:mu} and \eqref{eq:sigma}, respectively. The time span $\tau_{L}$ for the $L$ reactions follows the gamma distribution,  $\tau_{L} \sim \Gamma(L, 1/a_{0}(\vec{x}))$. The number of firings $k_{j}^{\mathcal{B}}$ for the reaction channel $R_{j}$, fired within the time span $\tau_{L}$, is sampled from a sequence of correlated binomial distributions,
%--------------------------------------------------------
\begin{equation}\label{eq:R:sampling}
k_{j}^{\mathcal{B}} \sim \mathcal{B}\left( L - \sum_{m=1}^{j-1} k_{m}^{\mathcal{B}}, \frac{ a_{j} (\vec{x})}{  a_{0}(\vec{x}) - \sum_{m=1}^{j-1} a_{m}(\vec{x})  }  \right)  \PERIOD
\end{equation}
%--------------------------------------------------------
This approach requires at most $M-1$ drawings of random numbers since $k_{M}^{\mathcal{B}} = L - \sum_{j=1}^{M-1} k_{j}^{\mathcal{B}}$. Furthermore, it can be shown that the sampling procedure is invariant under the permutation of reaction channels indices \cite{Auger:2006}. This fact can be exploited to reduce the number of samples drawn per time step by reordering the reactions indices in a way that the most probable reactions channels are sampled first. The $R$-leaping algorithm is summarised in Algorithm \ref{alg:R}.

The sampling of reaction channels from the bounded binomial distribution reduces the appearance of negative species, compared to sampling from the unbounded Poisson distribution. However, in systems involving species with population close to zero taking place in very fast reactions, the $R$-leaping algorithm might also introduce negative population. To control the appearance of negative population, an additional bounding condition for $L$ was proposed \cite{Auger:2006}. In systems with high rejection rates of the proposed state update, the total number of firings  is computed as $L = \min(L^\prime,L^{\prime\prime})$, where $L^\prime$ is given by \cref{eq:R:L} and
%--------------------------------------------------------
\begin{equation}\label{eq:R:L2}
L^{\prime\prime} = \min_{j=1,\ldots,M} \left( 1 - \theta \left(  1 - \frac{a_{0}(\vec{x})}{a_{j}(\vec{x})} \right) \right) L_{j}  \COMMA
\end{equation}
%--------------------------------------------------------
where
%--------------------------------------------------------
\begin{equation}
L_{j} = \min_{\substack{i=1,\ldots,N \\ \nu_{ij}<0}}  \left \lfloor \frac{ x_{i} }{ |\nu_{ij}| } \right \rceil \PERIOD
\end{equation}
%--------------------------------------------------------
The parameter $\theta$ controls appearance of negative species. Smaller values of $\theta$ lead to better control of negative species but also lead to  lower performance.
%
%
%
%
%~~~~~~~~~~~~~~~~~~~~~~~~~~~
\begin{algorithm}
\caption{\label{alg:R}R-Leaping}
{
\fontsize{8}{12}\selectfont
\begin{algorithmic}[1]
 \algsetup{indent=1.5em}
\STATE Initialise: $T_\textrm{end}$, $\vec{x} \gets \textbf{X}(0)$, $t \gets 0$, $steps\gets 0$, $p\gets$ frequency of reordering.
\WHILE{$t<T_\textrm{end}$}
	\STATE Compute $a_{j}(\vec{x})$ for $j=1, \dots , M$ and $a_{0}(\vec{x}) $%\gets \sum_{j=1}^{M} a_{j}( \vec{x})$
\IF{$\mod(steps,p)=0$}
\STATE Reorder the reactions such that $ a_{1}(\vec{x}) \geq a_{2}(\vec{x}) \geq \ldots \geq a_{M}(\vec{x}) $
\ENDIF
\STATE Compute $L$ by \cref{eq:R:L}, then set $L \gets \max(L, 1)$
\STATE Sample $k_{j}$ by \cref{eq:R:sampling} \label{alg:R:line:L}
\IF{ there is a negative component in $\vec{x} + \sum_{j=1}^{M} k_{j}\vec{\nu}_{j}$ }
\STATE  $L \gets L/2$ and go to \ref{alg:R:line:L}.
\ELSE{}
\STATE $\tau \gets \Gamma(L, 1/a_{0}(\vec{x}))$
\STATE $\vec{x} \gets \vec{x} + \sum_{j=1}^{M} k_{j} \vec{\nu}_{j}$
\STATE $t \gets t+\tau$
\STATE $steps=steps+1$
\ENDIF
\ENDWHILE
\end{algorithmic}
}
\end{algorithm}
%~~~~~~~~~~~~~~~~~~~~~~~~~~~
%==============================

%==================================================
%==========	Section: S-Leaping
%==================================================
\section{S-leaping}\label{sec:S}
Here, we propose the $S$-leaping, an algorithm which combines the advantages of the $\tau$-leaping and $R$-leaping algorithms. The $S$-leaping couples the efficient time step selection of the $\tau$-leaping with the effective binomial sampling of the $R$-leaping algorithm. The coupling of the algorithms is achieved in the following way. First, the time-step $\tau$ is selected according to \cref{eq:ExpTauTimeStep}. Then the total number of firings $L$ that will take place in the time interval $[t,t+\tau)$ is estimated. Since in the $\tau$-leaping each reaction channel is independently sampled as $k_{j}\sim \mathcal{P}(a_{j} (\vec{x})\tau)$, the total amount of all firings $L$ follows the Poisson distribution with parameter $a_{0}(\vec{x})\tau$, i.e.,
%--------------------------------------------------------
\begin{equation}\label{eq:S:L}
L(t) \sim \mathcal{P}(a_{0}(\vec{x}) \tau) \COMMA
\end{equation}
%--------------------------------------------------------
for $\vec{x}=\vec{X}(t)$.
Knowing the number of reactions that will take place in $[t,t+\tau)$, the firings of each channel $k_{j}$ can be sampled from the binomial distribution given by \cref{eq:R:sampling}. If the sampled $L$ is zero, it means the system will advance to the time $t=t+\tau$ without any changes since no reaction will be fired.  In this case the system can be further advanced by setting $L=1$ and $\tau\sim\Gamma(1, 1/a_{0}(\vec{x}))$ and proceeding with the $S$-leaping algorithm. Notice that this is just one step of SSA since the $\Gamma$ distribution with parameters $1$ and $1/a_{0}(\vec{x})$ is equal to exponential distribution with parameter $1/a_{0}(\vec{x})$.  The $S$-leaping algorithm is summarised in \cref{alg:S}.

To control the appearance of the negative species, the S-leaping algorithm can inherit the control mechanism from the $\tau$-leaping or $R$-leaping. Here we bound $L$ similarly as in the $R$-leaping method. In systems with high rejections rates, the total amount of firings is computed as $L= \min(L^\prime, L^{\prime\prime})$, where $L^{\prime}$ is given by \cref{eq:S:L} and $L^{\prime\prime}$ by \cref{eq:R:L2}. If $L^{\prime\prime}$ was chosen, then the time step $\tau$ should be recomputed as $\tau \sim \Gamma(L, 1/a_{0}(\vec{x}))$.

Thanks to the coupling of the two algorithms, the $S$-leaping performs always as well as the $\tau$-leaping or $R$-leaping algorithm. In the non-stiff systems, the $S$-leaping benefits from the efficient time step selection and might outperform the $R$-leaping method. On the other hand,  in the case of big and stiff systems, the $S$-leaping surpasses the $\tau$-leaping due to the effective sampling of the reaction channels. Moreover, since the behaviour of the system might change over time, the $S$-leaping can outperform both the $R$-leaping and $\tau$-leaping. Finally, since the $S$-leaping uses the same time-step selection as the $\tau$-leaping, the algorithm can easily be extended to an adaptive explicit-implicit version.

%~~~~~~~~~~~~~~~~~~~~~~~~~~~
\begin{algorithm}
\caption{\label{alg:S}S-Leaping}
{
\fontsize{8}{12}\selectfont
\begin{algorithmic}[1]
\STATE Initialise: $T_\textrm{end}$, $\vec{x} \gets \textbf{X}(0)$, $t \gets 0$, $steps\gets 0$, $p\gets$ frequency of reordering.
\WHILE{$t<T_\textrm{end}$}
	\STATE Compute $a_{j}(\vec{x})$ for $j=1, \dots , M$ and $a_{0}(\vec{x}) $%\gets \sum_{j=1}^{M} a_{j}( \vec{x})$
\IF{$\mod(steps,p)=0$}
\STATE Reorder the reactions such that $ a_{1}(\vec{x}) \geq a_{2}(\vec{x}) \geq \ldots \geq a_{M}(\vec{x}) $
\ENDIF
\STATE Compute $\tau$ by \cref{eq:ExpTauTimeStep}
\STATE Sample $L$ by \cref{eq:S:L}  \label{alg:s:line:tau}
\IF{L=0}
\STATE $t\gets t+\tau$
\STATE $L\gets 1$ and $\tau\gets \Gamma(1, 1/a_{0}(\vec{x})$
\ENDIF
\STATE Sample $k_{j}$ by \cref{eq:R:sampling}
\IF{ there is a negative component in $\vec{x} + \sum_{j=1}^{M} k_{j}\vec{\nu}_{j}$ }
\STATE  $\tau \gets \tau/2$ and go to \ref{alg:s:line:tau}
\ELSE{}
\STATE Update: $\vec{x} \gets \vec{x} + \sum_{j=1}^{M} k_{j}\vec{\nu}_{j}$
\STATE  $t \gets t+\tau$
\STATE $steps=steps+1$.
\ENDIF
\ENDWHILE
\end{algorithmic}
}
\end{algorithm}
%~~~~~~~~~~~~~~~~~~~~~~~~~~~
%
%
%
%
%
%
%==============================
% SubSec: Adaptive S
%==============================
\subsection{Adaptive S-Leaping}\label{sec:AdaptiveS}
The adaptive leap methods switch between explicit (\cref{alg:S}) and implicit method depending on the stiffness of the system. The implicit $S$-leaping method updates the system state as,
%-----------------------------------
\begin{eqnarray}\label{ImpSupdate}
\vec{X}(t+\tau) = \vec{X}(t) + \sum_{j=1}^{M} \vec{\nu}_{j} \, k_{j}^{\mathcal{B}}\left( \vec{X}(t+\tau)\right) \PERIOD
\end{eqnarray}
%-----------------------------------
This requires sampling random numbers $k_{j}^{\mathcal{B}}(\vec{X}(t+\tau))$ from the binomial distribution $\mathcal{B}(\alpha(\vec{x}^\prime),\beta(\vec{x}^\prime ))$ with mean and variance evaluated at the unknown state $\vec{x}^\prime=\vec{X}(t+\tau)$ given by,
%-----------------------------------
\begin{eqnarray}
\alpha(\vec{x}^\prime) &=&  L(\vec{x}^\prime) - \sum_{m=1}^{j-1} k^{\mathcal{B}}_{m}(\vec{x}^\prime) \COMMA \nonumber \\
\beta(\vec{x}^\prime) &=& \frac{  a_{j}(\vec{x}^\prime) }{ a_{0}(\vec{x}^\prime)  - \sum_{m=1}^{j-1} a_{m}(\vec{x}^\prime) } \PERIOD  \nonumber
\end{eqnarray}
%-----------------------------------
In the implicit $\tau$-leaping algorithm, each firing $k_{j}^{\mathcal{P}}$ is independently approximated by the partially implicit formulation given by \cref{eq:PoissSplit}. This can not be directly applied in the $S$-leaping, since each sample $k_{j}^{\mathcal{B}}$  depends on all previously drawn samples $k^{\mathcal{B}}_{\ell}, \,\, \ell =1,2, \ldots, j-1$.  
The partially implicit treatment for the $S$-leaping can be obtained by rather considering the distribution of the whole vector of all firings $(k_1,\ldots,k_M)$, i.e., the multinomial distribution with parameter \REVa{$\big(\frac{a_{1}(\vec{x})}{a_{0}(\vec{x})},\ldots, \frac{a_{M}(\vec{x})}{a_{0}(\vec{x})} \big )$ and $L$ the number of trials}.
\REVa{
If $k^{\mathcal{M}}_{j}$ is the $j$-th entry of a random vector that follows the multinomial distribution with parameters $\big (\frac{a_{1}(\vec{x})}{a_{0}(\vec{x})},\ldots, \frac{a_{M}(\vec{x})}{a_{0}(\vec{x})} \big)$ and $L$, 
}
then $k^{\mathcal{M}}_{j}$ can be expressed as the sum of a random variable with mean  $ \frac{a_{j}(\vec{x})}{a_{0}(\vec{x})} L$ and the zero mean variable  $k^{\mathcal{M}}_{j}  -  \frac{a_{j}(\vec{x})}{a_{0}(\vec{x})} L$, i.e.,
%-----------------------------------
\begin{equation}
k^{\mathcal{M}}_{j}=  \frac{a_{j}(\vec{x})}{a_{0}(\vec{x})} L + k^{\mathcal{M}}_{j} - \frac{a_{j}(\vec{x})}{a_{0}(\vec{x})} L \PERIOD
\end{equation}
%-----------------------------------
The variable  $\frac{a_{j} }{a_{0}} L$ is evaluated at the unknown state $\vec{X}(t+\tau)$, while the variable $ k^{\mathcal{M}}_{j} - \frac{a_{j}}{a_{0}}L $ is evaluated at the known state $\vec{X}(t)$. The partial implicit approximation to the variables $k_{j}^{\mathcal{M}}$ leads to the following system of implicit equations,
%-----------------------------------
\begin{equation}\label{eq:S:ImpSsystem}
\vec{x}^\prime = \vec{x} + \sum_{j=1}^{M} \vec{\nu}_{j} \frac{a_{j}(\vec{x}^\prime)}{a_{0}(\vec{x}^\prime)}L(t+\tau)  + \sum_{j=1}^{M} \vec{\nu}_{j} \left( k_{j}^{\mathcal{M}}(\vec{x})  - \frac{a_{j}(\vec{x})}{ a_{0}(\vec{x}) }L(t) \right) \COMMA
\end{equation}
%-----------------------------------
for $\vec{x}=\vec{X}(t)$ and $\vec{x}^\prime=\vec{X}(t+\tau)$. %
%............................
\REVb{
Since the multinomial random vectors $k_j^\mathcal{M}$ in \cref{eq:S:ImpSsystem} depend on the known state $\vec{x}$ and since the $j$-th element of the multinomial distribution follows binomial distribution, $k_j^{\mathcal{M}}(\vec{x})$ are computed by \cref{eq:R:sampling}.
}
%............................
However, $L(t+\tau)$ is also a random variable from Poisson distribution with the parameter evaluated at the unknown state $\vec{X}(t+\tau)$,
%-----------------------------------
\begin{equation}
L(t+\tau) \sim \mathcal{P}(a_{0}(\vec{X}(t+\tau))\tau) \PERIOD
\end{equation}
%-----------------------------------
The term $L(t+\tau)$ could be expressed in the partial implicit manner following \cref{eq:PoissSplit}. However, a simpler formulation can be obtain by a mean approximation,
%-----------------------------------
\begin{equation}\label{eq:AS:L}
L(t+\tau) \approx a_{0}(\vec{X}(t+\tau)) \tau \COMMA
\end{equation}
%-----------------------------------
which corresponds to the computation of the $L$ in the $R$-leaping method. %
%........................
\REVb
{
The advantage of the mean approximation in \cref{eq:AS:L} is that it significantly reduces numerical complexity of the implicit system in \cref{eq:S:ImpSsystem}, while the partial implicit approximation provided by \cref{eq:PoissSplit} would increase the complexity even more. Since we are dealing with stiff system, increased complexity could reduce accuracy of the numerical methods used for solving the implicit system of equations.
}
%........................
If $\vec{X}^{\star}$ is the solution of the implicit system of Eq. \eqref{eq:S:ImpSsystem}, then the implicit update is obtained as
%-----------------------------------
\begin{equation}
\vec{X}(t+\tau) = \vec{X}(t) + \sum_{j=1}^{M} \vec{\nu}_{j} k_{j}^{\mathcal{M}\star} \COMMA
\end{equation}
%-----------------------------------
where
%-----------------------------------
\begin{equation}\label{eq:as:kj}
k_{j}^{\mathcal{M}\star} =  \left \lfloor  a_{j}(\vec{X}^{\star}) \tau + k_{j}^{\mathcal{M}}(\vec{X}(t)) - \frac{a_{j}(\vec{X}(t))}{a_{0}(\vec{X}(t))} L(t)  \right \rceil  \PERIOD
\end{equation}
%-----------------------------------
This means that both, the implicit $\tau$-leaping and implicit $S$-leaping algorithm solve the implicit system with the same computational complexity. 
%
%...................
\REVb{
However, since the implicit $S$-leaping can exploit reordering of reaction channels, it might result in less random number generations (at most $M$ samples) than the implicit $\tau$-leaping (always $M$ samples). This might allow the implicit $S$-leaping to outperform the implicit $\tau$-leaping, especially in large stiff systems where only few reaction channels are fired per time step.
}
%...................
The adaptive $S$-leaping method is summarised in \cref{alg:AdaptS}.
%
%
%~~~~~~~~~~~~~~~~~~~~~~~~~~~
\begin{algorithm}
\caption{\label{alg:AdaptS} Adaptive S-Leaping}
{
\fontsize{8}{12}\selectfont
\begin{algorithmic}[1]
\STATE Initialise: $T_\textrm{end}$, $\vec{x} \gets \textbf{X}(0)$, $t \gets 0$, $steps\gets 0$, $p\gets$ frequency of reordering.
\WHILE{$t<T_\textrm{end}$}
	\STATE Compute $a_{j}(\vec{x})$ for $j=1, \dots , M$ and $a_{0}(\vec{x})$% \gets \sum_{j=1}^{M} a_{j}( \vec{x})$
\IF{$\mod(steps,p)=0$}
\STATE Reorder the reactions such that $ a_{1}(\vec{x}) \geq a_{2}(\vec{x}) \geq \ldots \geq a_{M}(\vec{x}) $
\ENDIF
\STATE Compute $\tau^{(\textrm{ex})}$ by \cref{eq:ExpTauTimeStep} and $\tau^{(\textrm{im})}$ by \cref{eq:ImplTau}
\IF{$ \tau^{(\textrm{im})} > 100\; \tau^{(\textrm{ex})} $}
		\STATE System is stiff and $\tau \gets \tau^{(\textrm{im})}$
	\ELSE{}
		\STATE System is non-stiff and $\tau \gets \tau^{(\textrm{ex})}$
	\ENDIF
    \IF{ the system is currently non stiff} \label{alg:as:line:tau}
    \STATE Compute $L$ by \cref{eq:S:L}
     \IF{L=0}
	\STATE $t\gets t+\tau$
	\STATE $L\gets 1$ and $\tau\gets \Gamma(1, 1/a_{0}(\vec{x})$
	\ENDIF
	\STATE Sample  $k_{j}$ by \cref{eq:R:sampling} 
    \ELSE{}
    \STATE Compute $k_{j}$ by \cref{eq:as:kj}, where $L(t)$ is given by \cref{eq:S:L}
    \ENDIF
	\IF{ there is a negative component in $\vec{x} + \sum_{j=1}^{M} k_{j} \vec{\nu}_{j}$  }
	\STATE  $\tau \gets \tau/2$ and go to \ref{alg:as:line:tau}
	\ELSE{}
	\STATE $ \vec{x} \gets \vec{x} + \sum_{j=1}^{M} k_{j} \vec{\nu}_{j}$
	\STATE $t \gets t+\tau$
	\STATE $steps=steps+1$.
	\ENDIF
\ENDWHILE
\end{algorithmic}
}
\end{algorithm}
%~~~~~~~~~~~~~~~~~~~~~~~~~~~
%==================================================
%
%
%
%
%

%==================================================
%========	 Section: Numerical Simulations   	     ============
%==================================================

\section{Numerical Simulations}\label{sec:numerics1}
To demonstrate the efficiency of the $S$-leaping algorithm, it is studied in comparison with the $\tau$-leaping and $R$-leaping methods on four reaction networks. The first one is a non-stiff system simulating decaying dimerization. The second system is a stiff decaying dimerization with reversible reaction channels in partial equilibrium. The third one is a system with very fast dynamics simulating the evolution of Bacillus subtilis. The last one is a LacZ/LacY system, which consists of a relatively large amount of reactions and which stiffness change over time. For each reaction network and each algorithm we measure two quantities: the error and the execution time of the algorithm. 

The error is measured as follows. For 25 equally distributed time points and all species we estimate the distance between the distributions of the tested algorithm and the SSA \cite{Cao:2006}. 
\REVb{Since the methods do not advance the system with a fixed time step, the population at a given time is approximated by the population at the closest time where the algorithm has landed.}
The distance $d$ between two distributions $P$ and $Q$ is approximated by the estimated histogram as
\begin{eqnarray}
d=\Delta \sum_{k}  | \tilde P(k) - \tilde Q(k) | \COMMA
\end{eqnarray}
where $\Delta$ is the bin size and $\tilde P(k),\tilde Q(k)$ are the values of the histogram for $P$ and $Q$ at the $k$-th bin. The histogram is computed using $N_s = 10^4$ independent trajectories and $K=10$ number of bins. Finally, the average error over all time points and all species is reported. 
%..........................
\REVb{
This definition of the error can be interpreted as a global error, since it accounts for temporal and interspecies error of the algorithm. Averaging the error over many time points takes into account the error not only at equilibrium but at transient regimes as well.
}

\REVb{
In \cite{Cao:2006b} the histogram self-distance was introduced as a measure of accuracy of the histogram distance estimation. Any estimate bellow the value of self-distance should be considered inaccurate. A bound for the self-distance was derived in \cite{Cao:2006b} and is given by $\sqrt{{4K}/(\pi N_s)}$. In all the histogram error plots we show the self-distance estimate as a constant blue line. Although the errors close or bellow this line should not be considered accurate we present them for completeness. 
}

\REVb{
The execution time is averaged over $10$ independent runs. The ratio between the execution time of SSA and the execution time of each algorithm is reported as a speed-up. 
Note that here the SSA is  used only as a reference in order to compare the relative speed-up of the three approximate algorithm, the $\tau$-, the $R$- and $S$-leaping. The execution times of the three methods are compared using as a reference the execution time of SSA. Hence changing the base implementation of SSA will not affect these comparisons results. 
Moreover, since we report the speed-up over the baseline SSA rather than CPU time, the presented results do not depend on the type of the used computer.
}

All reaction networks discussed in this section follow the law of mass action which states that the rate of a reaction is proportional to the product of the concentrations of the reactants. For example, for the reaction $2S_1 + S_2 \rightarrow S_3$ with reaction rate $c$, the propensity  is defined as,
\begin{equation}
a(\vec{x}) = c\, x_1 \, (x_1-1) \, x_2 \PERIOD
\end{equation}
For the general formula of the propensity function under the law of mass action we refer to \cite{Anderson:2011,Erban:2007}.

\REVb{
All the tested methods are implemented in the {\CC}  language, using the {\CC}11 random number generator library and the code is publicly available\footnote{\url{https://github.com/JanaLipkova/SSM}}.
}

%~~~~~~~~~~~~~~~~~~~~~~~~~~~
\begin{table}[t]
\begin{center}
\renewcommand*{\arraystretch}{1.1}
\begin{tabular}{ c | c | c  | c }
& Reaction   & Reaction Rate   &  Reaction Rate  \\
&                 &  (non-stiff)           &    (stiff)               \\ \hline
 $R_{1}$  &   $ S_{1} $   $ {\longrightarrow} $   $  \emptyset $      &  1           &  1   \\  \hline
 $R_{2}$  &    $ S_{1}+S_{1}  $   $ {\longrightarrow} $   $  S_{2} $  &  0.002    &  10  \\  \hline
 $R_{3}$  &    $ S_{2}  $   $ {\longrightarrow} $   $  S_{1}+S_{1} $  &  0.5        &   1000 \\  \hline
 $R_{4}$  &    $ S_{2} $   $  {\longrightarrow} $   $ S_{3} $              &  0.04      &  0.1
\end{tabular}
\end{center}
\caption{The reaction network for the Dimerization system studied in \cref{sec:nonstiff:dim,sec:stiff:dim}.}
\label{table:dim}
\end{table}
%~~~~~~~~~~~~~~~~~~~~~~~~~~~

%==============================
%=== Subsection: Test 1
%==============================
\subsection{Non-stiff Decaying Dimerization}\label{sec:nonstiff:dim}
\begin{figure}
\centering
  	\includegraphics[width=0.32\textwidth]{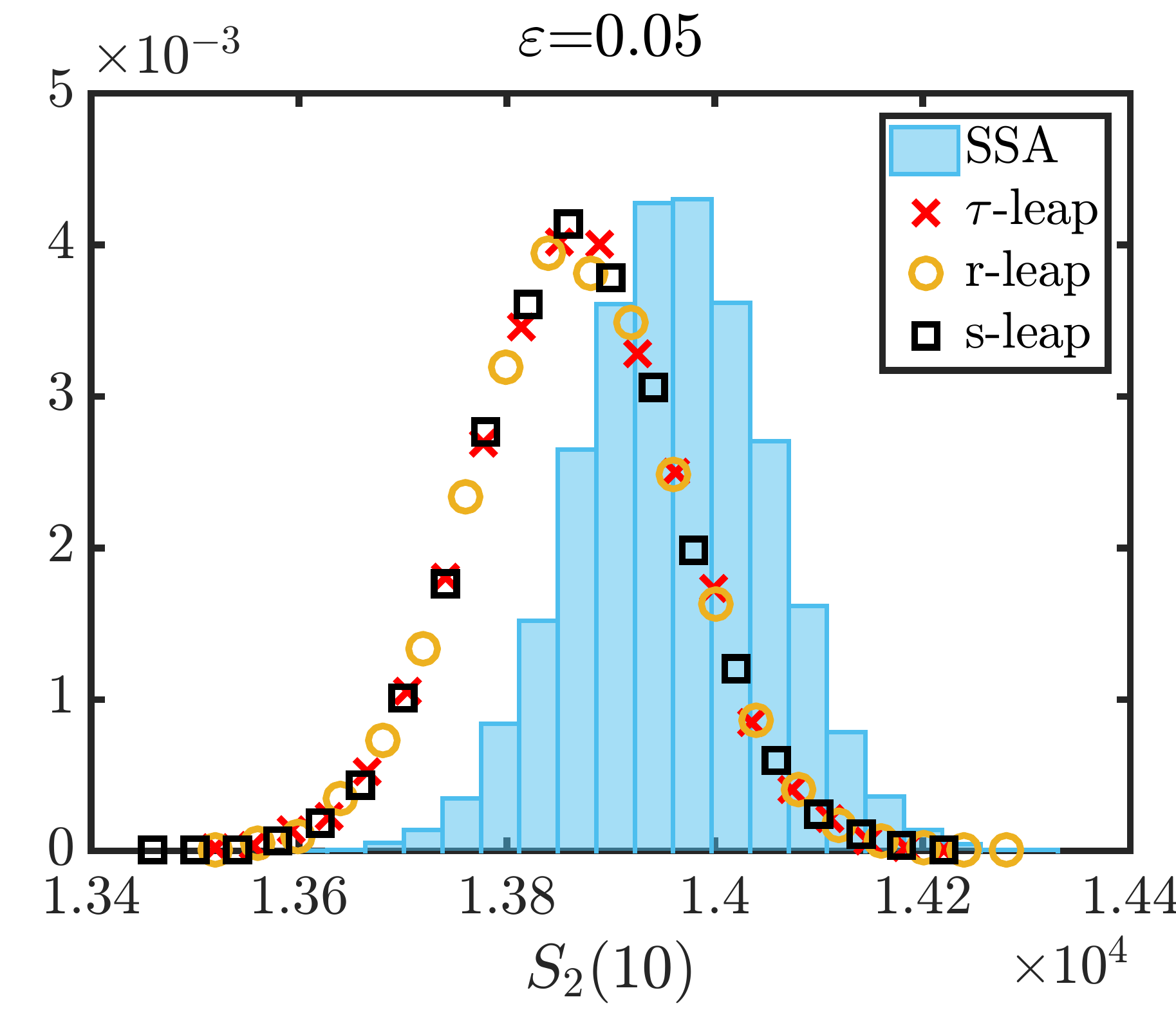}
	\includegraphics[width=0.32\textwidth]{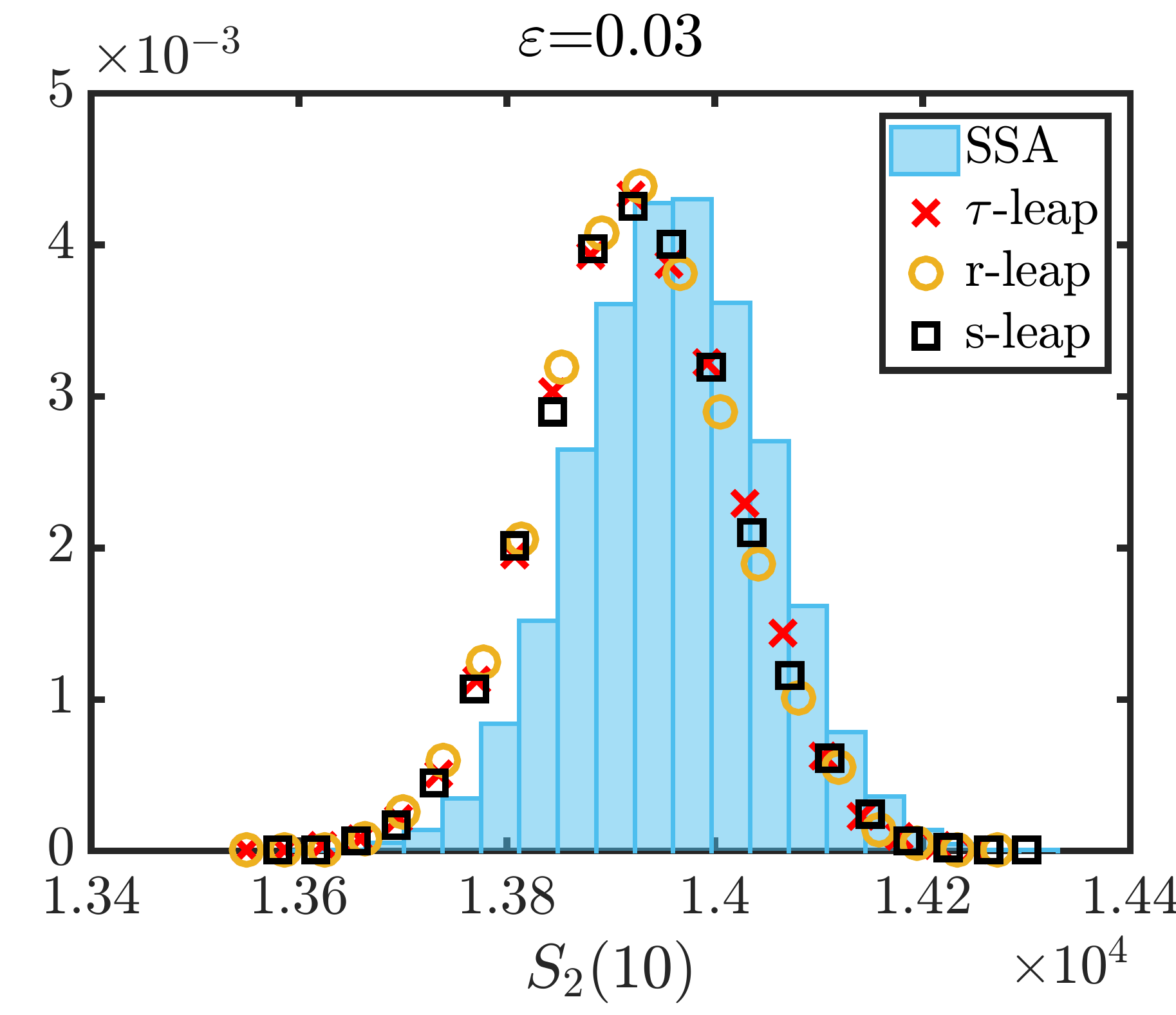}
   	\includegraphics[width=0.32\textwidth]{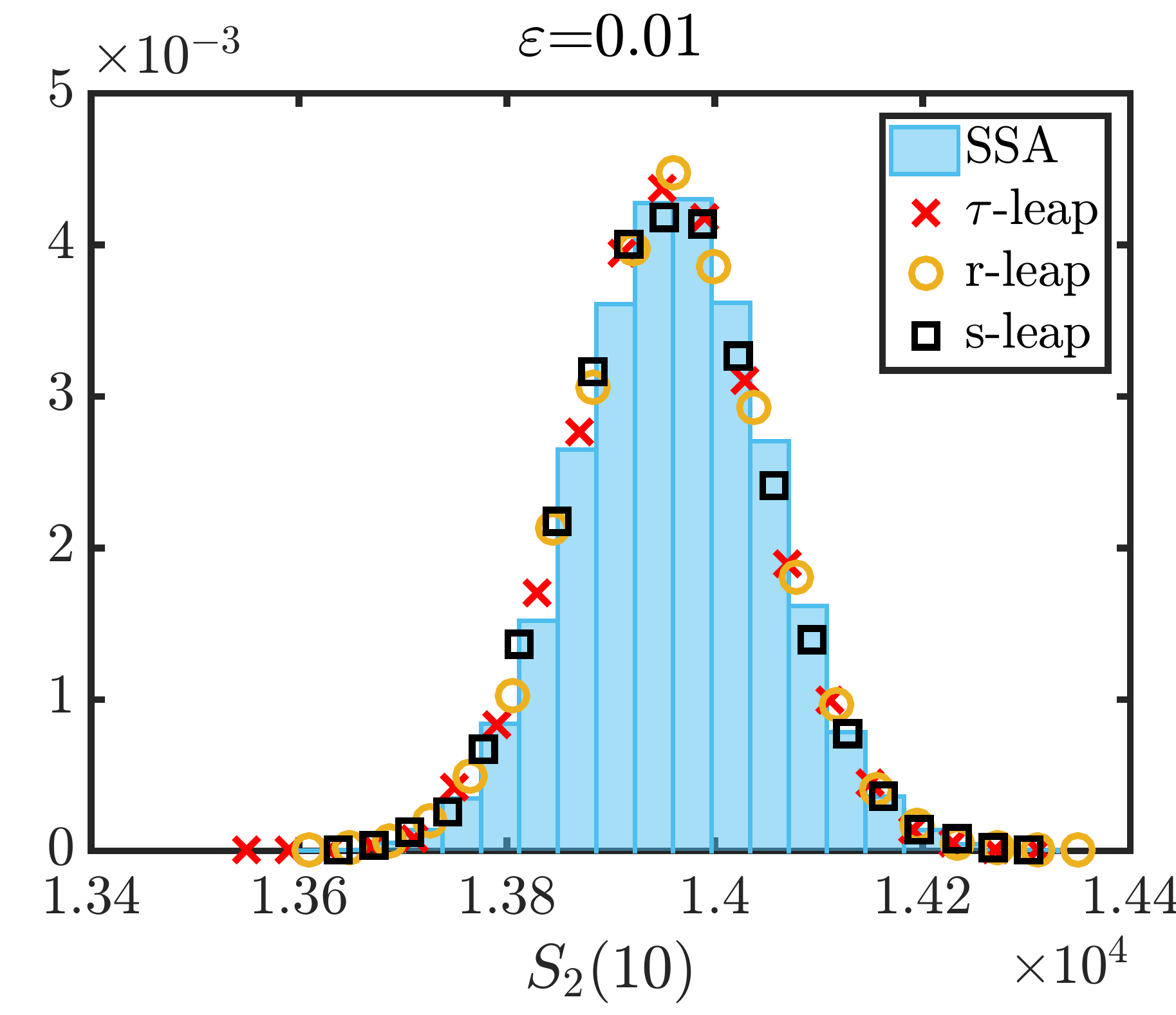}
\caption{Convergence of the approximate leap solutions to the exact SSA solution with decreasing values of the accuracy parameter $\varepsilon$ for the non-stiff dimerization system of \cref{sec:nonstiff:dim}.}
\label{fig:dim:conv}
\end{figure}

\begin{figure}
\centering
  	\includegraphics[width=0.45\textwidth]{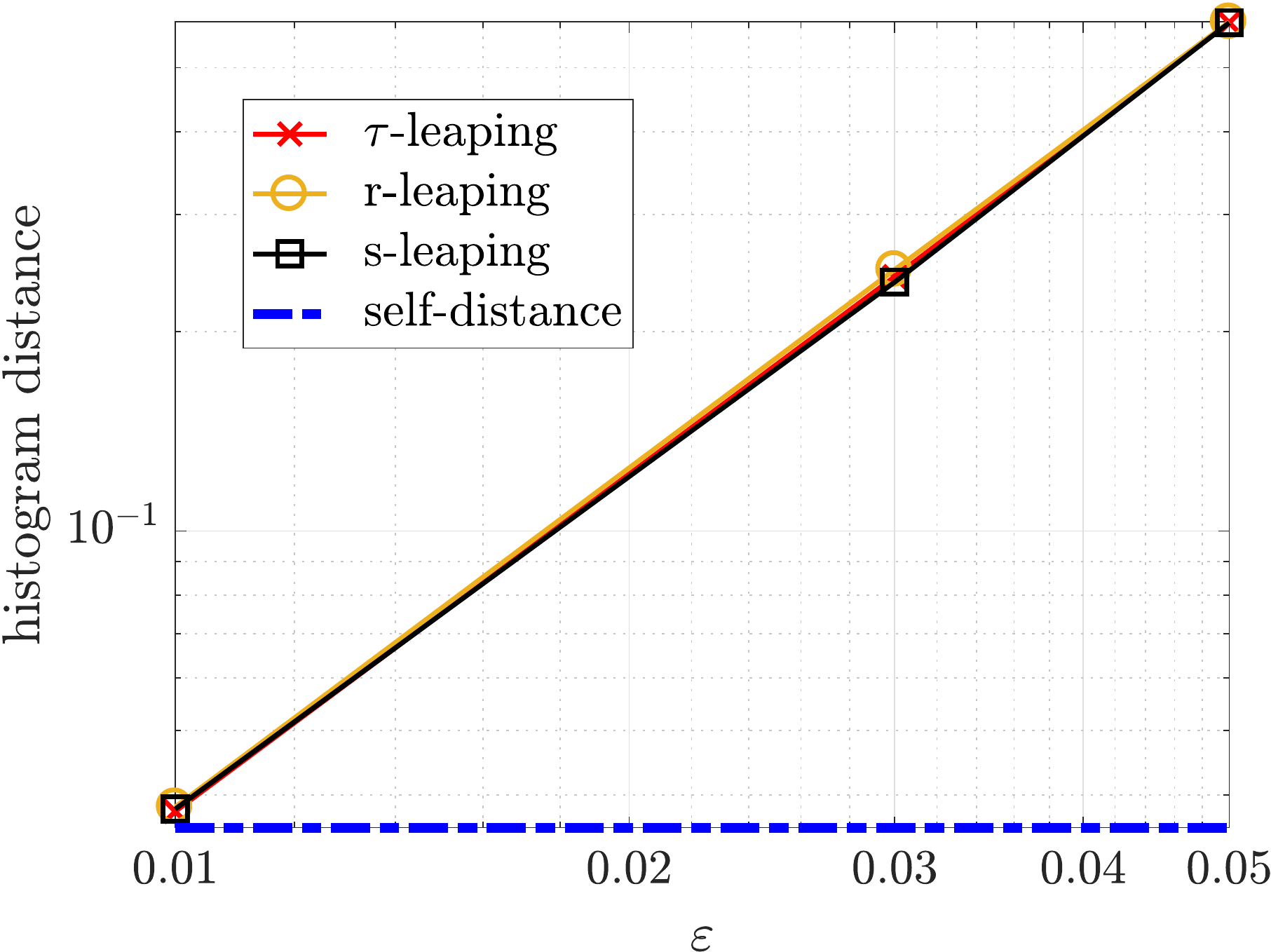}
	\includegraphics[width=0.45\textwidth]{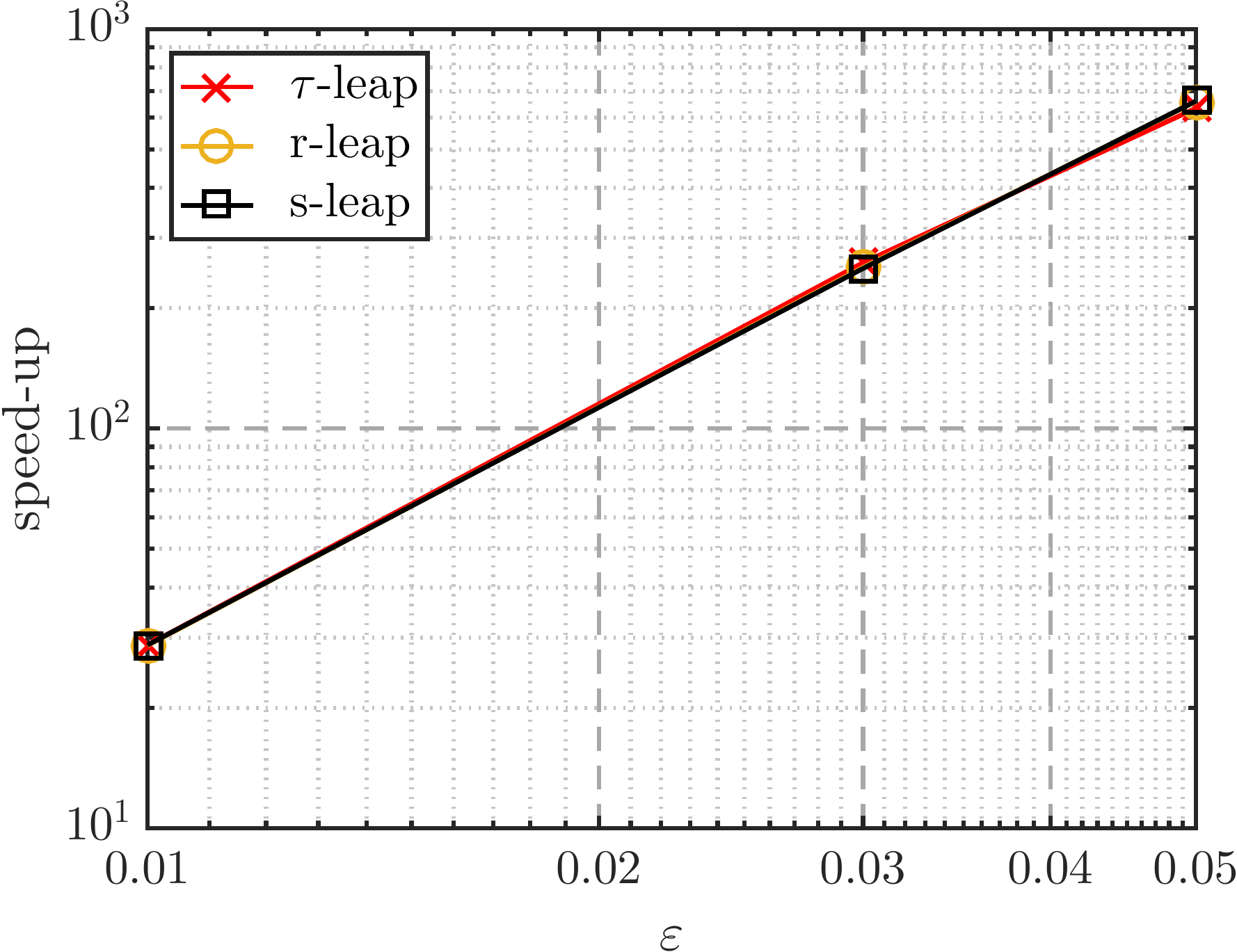}
\caption{Errors and efficiency for non-stiff dimerization system of \cref{sec:nonstiff:dim}.}
\label{fig:dim:ns:error:time}
\end{figure}

Following the same test as in \cite{Auger:2006,Gillespie:2001,Gillespie:2003} we consider the non-stiff decaying dimerization system summarised in \cref{table:dim}. The initial populations are $\vec{X}(0)=( 4150, 39565,3445)$ and the system is evolved until $T_{\textrm{end}}=10$ using $\varepsilon=\{0.01, 0.03,0.05\}$. In \cref{fig:dim:conv} we show the convergence of the histograms of the approximate algorithms to that of the SSA for the second species $S_{2}$   at time $t=10$. All the approximate methods converge to the SSA solution as the accuracy parameter $\varepsilon$ decreases. \cref{fig:dim:ns:error:time} shows the accuracy (left) and the performance (right) for all leap methods with different accuracy parameter $\varepsilon$. In this system, all the leaping methods have comparable accuracy and performance. No additional speed-up was obtained by reordering of the reaction channels in the $R$-leaping and $S$-leaping, since in each step of the simulation all reaction channels are fired.

%~~~~~~~~~~~~~~~~~~~~~~~~~~~
%~~~~~~~~~~~~~~~~~~~~~~~~~~~

%~~~~~~~~~~~~~~~~~~~~~~~~~~~
%~~~~~~~~~~~~~~~~~~~~~~~~~~~

%==============================
%=== Subsection: Stiff Dimerization
%==============================
\subsection{Stiff  Decaying Dimerization } \label{sec:stiff:dim}

\begin{figure}
\centering
  	\includegraphics[width=0.45\textwidth]{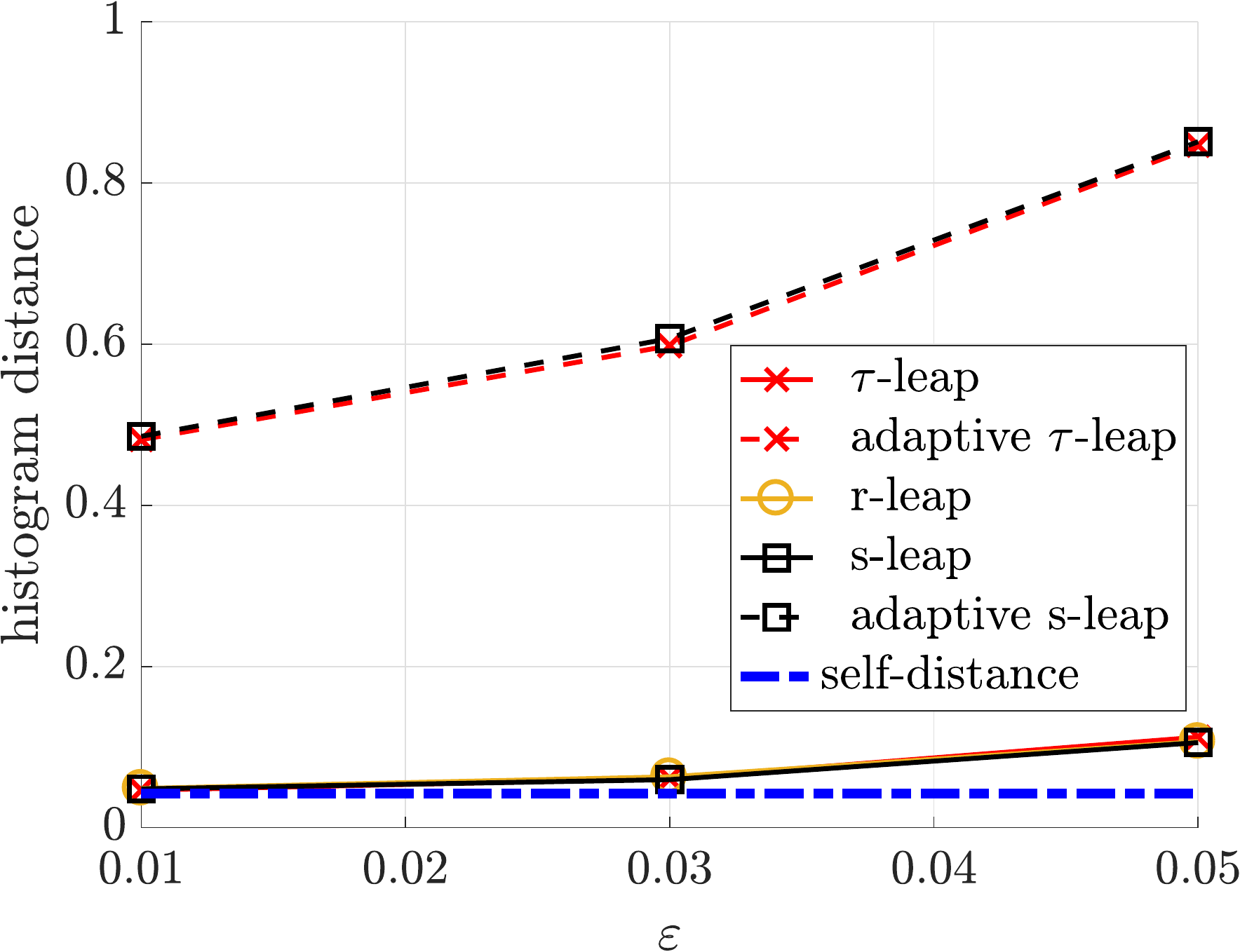}
	\includegraphics[width=0.45\textwidth]{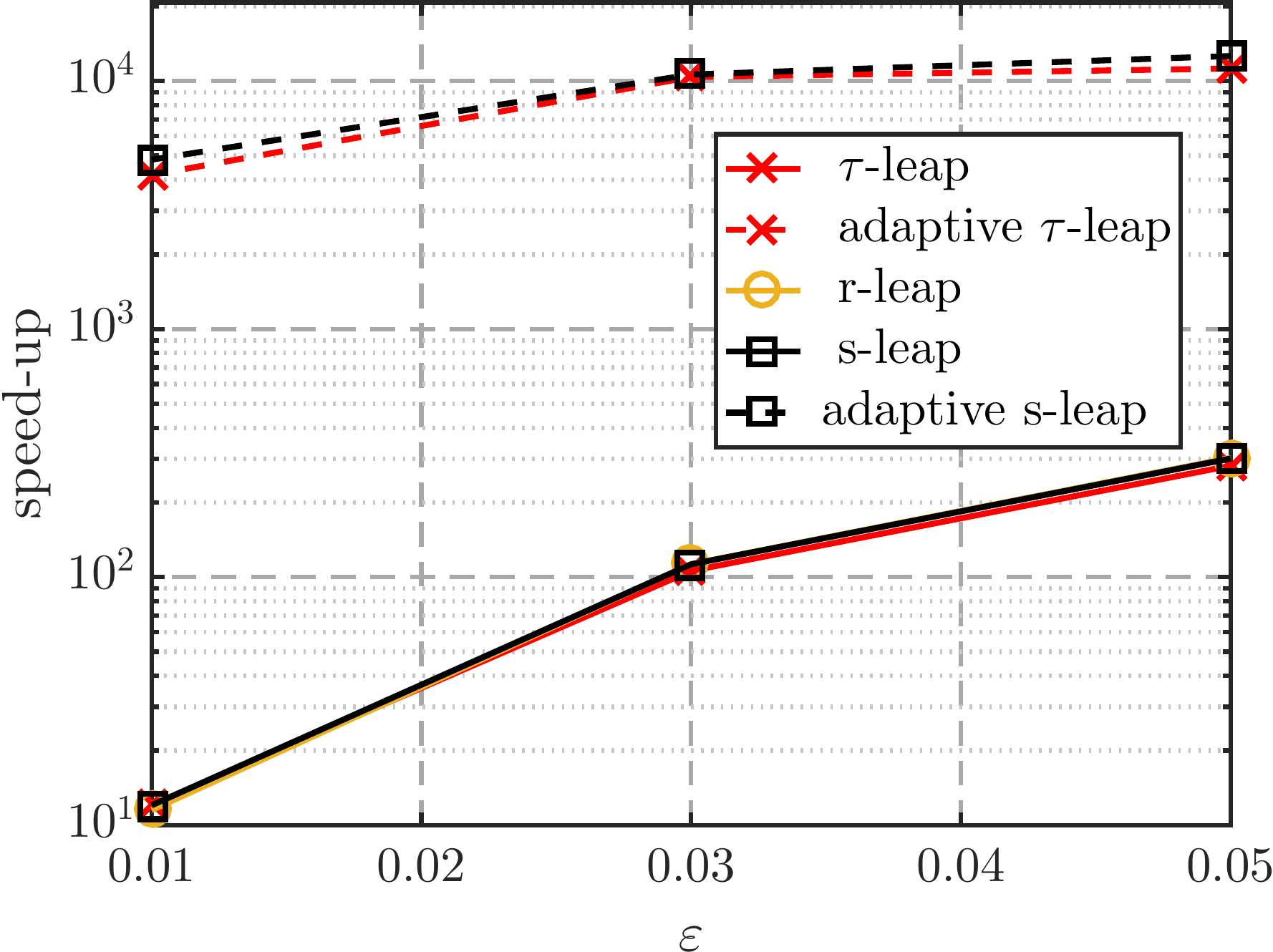}
\caption{Errors and efficiency for stiff dimerization system discussed in \cref{sec:stiff:dim}. }
\label{fig:dim:s:error:time}
\end{figure}

To study the efficiency of the adaptive $S$-leaping method we consider the stiff decaying dimerization system studied in \cite{Cao:2007,Rathinam:2003}. The system is defined by the same set of reactions and initial conditions as in Section \ref{sec:nonstiff:dim}, see Table \ref{table:dim}. The stiffness arises from the reaction rates that vary by a few orders of magnitude. The behaviour of this system changes over time starting with a non-stiff phase. However, once the reversible reactions $R_2$ and $R_3$ approach the equilibrium, the system becomes stiff. 
%.........................................
\REVb{
Under this set up $S_1$ and $S_2$ are the fast variables, while $S_3$ is the slow variable.
}
%......................................... 
The system is evolved until the final time $T_{\textrm{end}}=10$  for $\varepsilon = \{0.01, 0.03,0.05\}$.
\\
In \cref{fig:dim:s:error:time} we present the accuracy and the performance of the adaptive $\tau$-leaping and adaptive $S$-leaping as well as the explicit $R$-leaping, $\tau$-leaping and $S$-leaping. All explicit methods reach comparable accuracy and performance.
\REVb{
The adaptive methods provide significant speed-up over their explicit counterparts. The reduced accuracy of the adaptive methods arise from the dumping effect of the implicit methods on the fast variables. As reported in \cite{Rathinam:2003}, the implicit schemes capture the distribution of the slow variable $S_3$ correctly. However, for the fast variables $S_1$ and $S_2$, the mean is computed correctly but the histogram distribution around the mean is too narrow. In \cite{Rathinam:2003} a downshifting strategy was proposed to restore the natural fluctuations in the fast variables by simulating the final time steps of the adaptive method with the explicit method. As shown in \cite{Cao:2005,Rathinam:2003}, the downshifting leads to correct histogram distributions for all variables at the final time. Since the downshifting procedure corrects the dumping effect only in the final time, the global error of the adaptive method will not be reduced to the level of the explicit methods. Since we report the global error, the downshifting strategy was not applied here. However, the downshifting procedure can be used to increase the accuracy of the adaptive $\tau$-leaping and adaptive $S$-leaping method at the fixed time point.
}
%.........................................
% The reduced accuracy in the adaptive method is compensated by a significant speed-up as reported in \cite{Cao:2005,Cao:2007}. Adaptive methods provide good alternative if fast approximate solutions are desired. For more accurate results the explicit methods should be used. All explicit methods reach comparable accuracy and performance. Reordering of the reaction channels did not provide additional speed-up, since all the reaction channels are being fired in each simulation step.
%
%
%
%~~~~~~~~~~~~~~~~~~~~~~~~~~~
%~~~~~~~~~~~~~~~~~~~~~~~~~~~
%
%
%==============================
%
%
%
%==============================
%=== Subsection: BSubtilis 
%==============================

\subsection{\textit{Bacillus subtilis}}\label{sec:bsub}

\begin{figure}
\centering
	\includegraphics[width=0.45\textwidth]{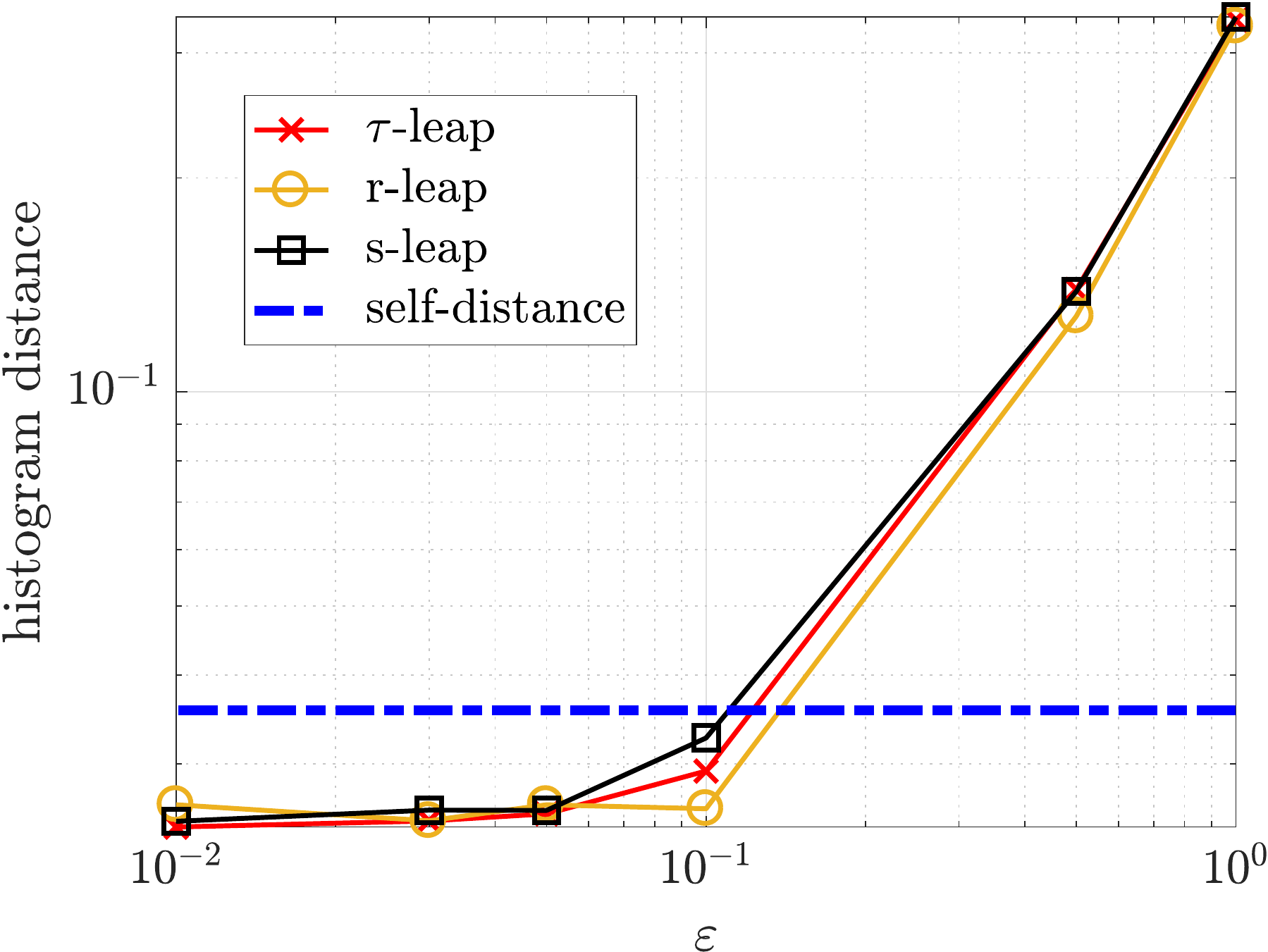}
  	\includegraphics[width=0.45\textwidth]{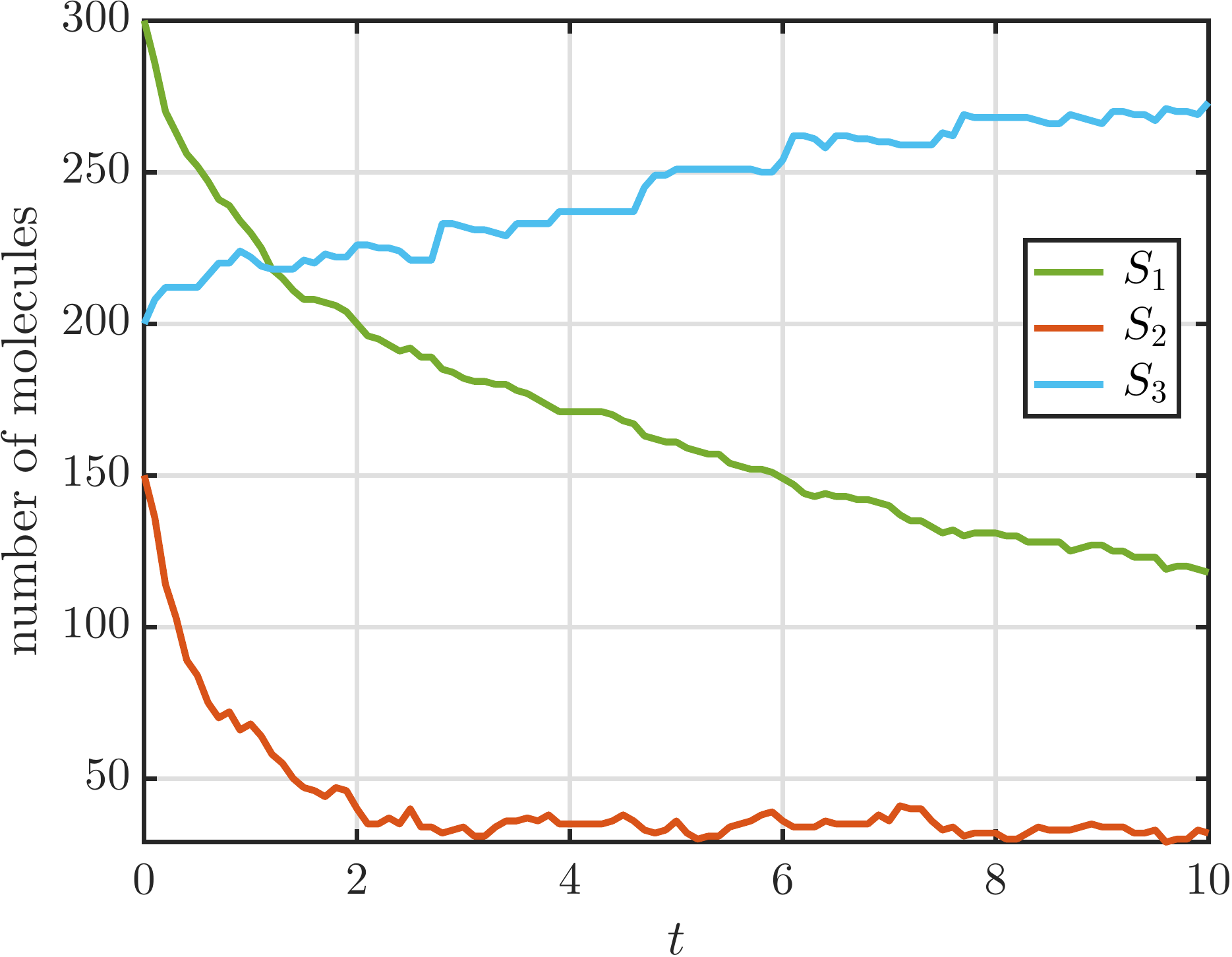}
\caption{Error of the leap methods (left) and a single trajectory of the Bacillus subtilis computed with SSA (right).}
\label{fig:bsub}
\end{figure}

This systems describes the cellular differentiation dynamics of the Bacillus subtilis which exhibits stochastic behaviour at the single-cell level \cite{Maamar:2007,Suel:2006}.  The differentiation dynamics depends on the expression of the transcriptional genes $S_{1}$=Spo0A, $S_{2}$=ComG and $S_{3}$=sinI and the reaction network is presented in \cref{table:bsub} \cite{Chattopadhyay:2013}. The system is evolved until the final time $T_{\textrm{end}}=10$ with initial population $\vec{X}(0)=(300,150,200)$. \cref{fig:bsub} (right) shows a single realisation of the Bacillus subtilis system computed with SSA. This system exhibits diverse reaction rates and very fast dynamics. As a consequence, the leap methods are strongly restricted by the leap condition and do not provide additional speed-up over SSA. The \cref{table:bsub:steps} shows the average number of steps executed by each method and the averaged CPU time for $\varepsilon=0.05$. The $R$-leaping algorithm advances the system only with one reaction per time step, emulating the SSA. Moreover, the $\tau$-leaping executes two times more steps than SSA. In this system, the $S$-leaping is the only method which requires less steps than SSA. Since in this case the leap methods do not provide additional speed up, the SSA alone would be the best choice. This example however shows that the $S$-leaping maintains its performance even in fast dynamical systems and outperforms the other leap methods. All methods reach comparable accuracy as shown in \cref{fig:bsub} (left).

%~~~~~~~~~~~~~~~~~~~~~~~~~~~
\begin{table}[b]
\begin{center}
\renewcommand*{\arraystretch}{1.15}
\begin{tabular}{ c | c | c  }
& Reaction   & Reaction Rate     \\ \hline
 $R_{1}$  &    $  \emptyset $        $ {\longrightarrow} $   $  S_{1}+3 S_{3} $  &  $ 1.51 \times 10^{-1}$    \\  \hline
 $R_{2}$  &    $ S_{1}+S_{2}  $    $ {\longrightarrow} $   $  4 S_{3}            $  &  $ 3.1 \times 10^{-4} $     \\  \hline
 $R_{3}$  &    $ S_{2}  $               $ {\longrightarrow} $   $  4 S_{3}            $  &  $ 3.4 \times 10^{-3} $     \\  \hline
 $R_{4}$  &    $ S_{3}  $               $ {\longrightarrow} $   $  S_{1}+S_{2}    $  &  $ 2.0 \times 10^{-2} $     \\  \hline
 $R_{5}$  &    $ S_{1}+2 S_{2} $  $ {\longrightarrow} $   $  \emptyset        $  &  $ 6.2 \times 10^{-5} $      \\  \hline
 $R_{6}$  &    $ 2 S_{1}           $   $ {\longrightarrow} $   $  S_{1}+S_{2}    $  &   $4.9 \times 10^{-4} $   
\end{tabular}
\end{center}
\caption{The reaction network for the Bacillus subtilis system studied in \cref{sec:bsub}.}
\label{table:bsub}
\end{table}
%~~~~~~~~~~~~~~~~~~~~~~~~~~~
%
%
%
%~~~~~~~~~~~~~~~~~~~~~~~~~~~
\begin{table}
\begin{center}
\renewcommand*{\arraystretch}{1.2}
\begin{tabular}{ c | c | c  }
Method          & Average number of steps   & Average CPU time [sec]    \\ \hline
SSA               &    266.6                               &  0.028                                  \\  \hline
$\tau$-leap    &    423.4                               &  0.086                                  \\  \hline
$R$-leap        &   263.2                                &  0.053                                  \\  \hline
$S$-leap        &   220.8                                &  0.045                                  \\ 
\end{tabular}
\end{center}
\caption{The averaged number of simulation steps and the execution time for the Bacillus subtilis system with the $\varepsilon=0.05$ presented in \cref{sec:bsub}.}
\label{table:bsub:steps}
\end{table}
%~~~~~~~~~~~~~~~~~~~~~~~~~~~
%
%
%
%~~~~~~~~~~~~~~~~~~~~~~~~~~~
%~~~~~~~~~~~~~~~~~~~~~~~~~~~
%
%==============================
%
%
%
%
%
%
%
%==============================
%=== Subsection: LacZ/LacY
%==============================
\subsection{LacZ/LacY}
\label{sec:lac}

In this section we consider the LacZ/LacY model which describes the expression of the LacZ and LacY genes and the activity of LacZ and LacY proteins in Escherichia Coli \cite{Kierzek:2002}. The reaction network consists of 22 reactions and 23 species. We present the reaction network, along with the reaction rate of each reaction in \cref{table:lac}. The propensity functions of this system vary by a few orders of magnitude making the system stiff. Moreover, the reaction system is considered inside a growing cell, with generation time $T_{\mathrm{gen}}=2100$. The growing volume changes the stiffness of the system over time since the propensities of the second and higher order reactions have to be rescaled by the volume. We consider two different initial conditions. In the first case we assume a small initial population where all species are initially 0 except  for PLac=1. In the second case, we consider bigger initial populations with all species initialized at 50 and PLac=100. In addition, the number of the species RNAP and ribosome are sampled every time step from a normal distribution $\mathcal{N}( 35(1 + t/T_{\mathrm{gen}} ), 3.5^2 )$ and $\mathcal{N}( 350(1 + t/T_{\mathrm{gen}} ), 35^2 )$, respectively for each case. The role of the system with small initial population is to investigate the behaviour of all methods in the presence of negative population, while the behaviour without the appearance of negative populations is studied in the system with the bigger initial population. 

\begin{figure}
\centering
  	\includegraphics[width=0.45\textwidth]{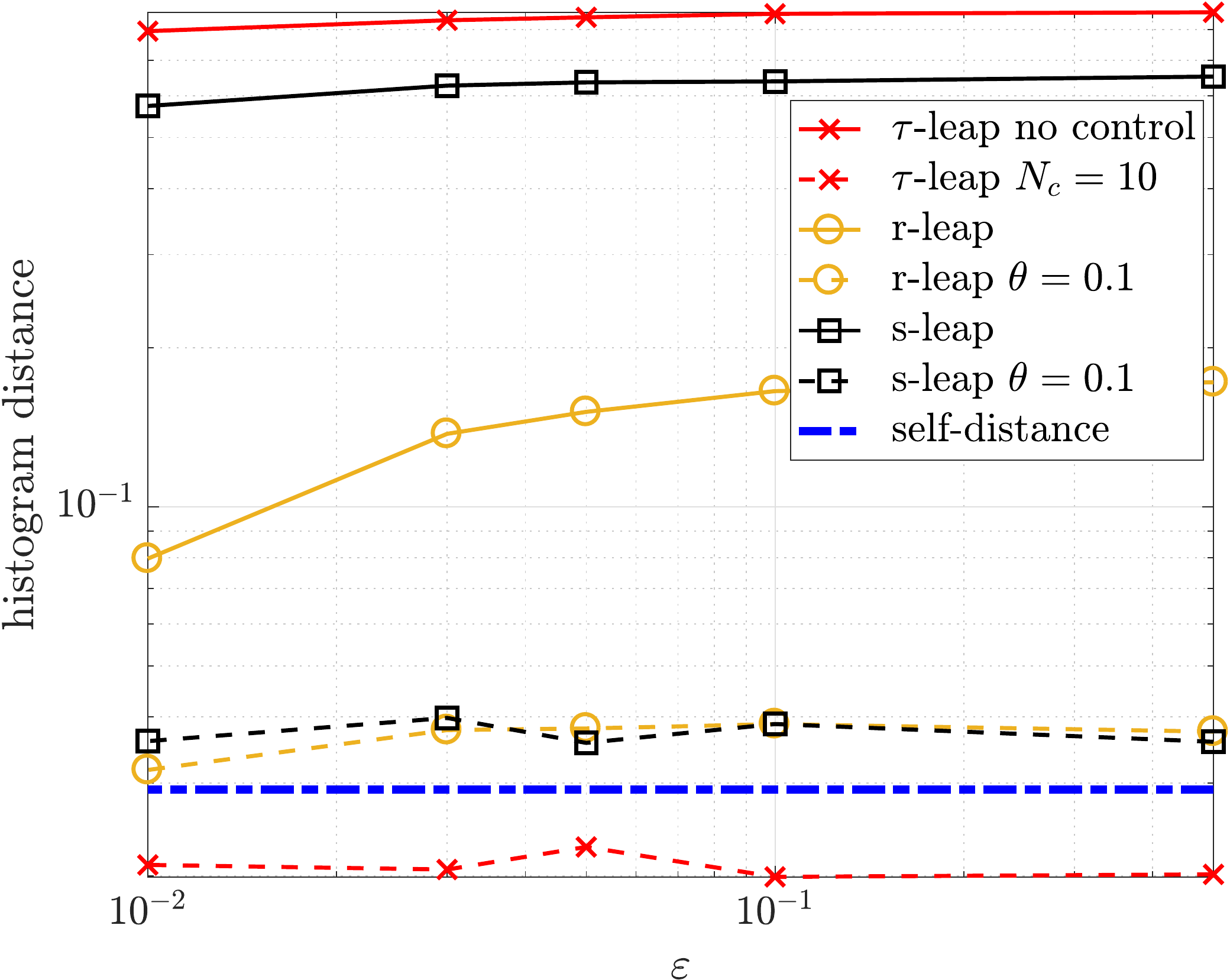}
	\includegraphics[width=0.45\textwidth]{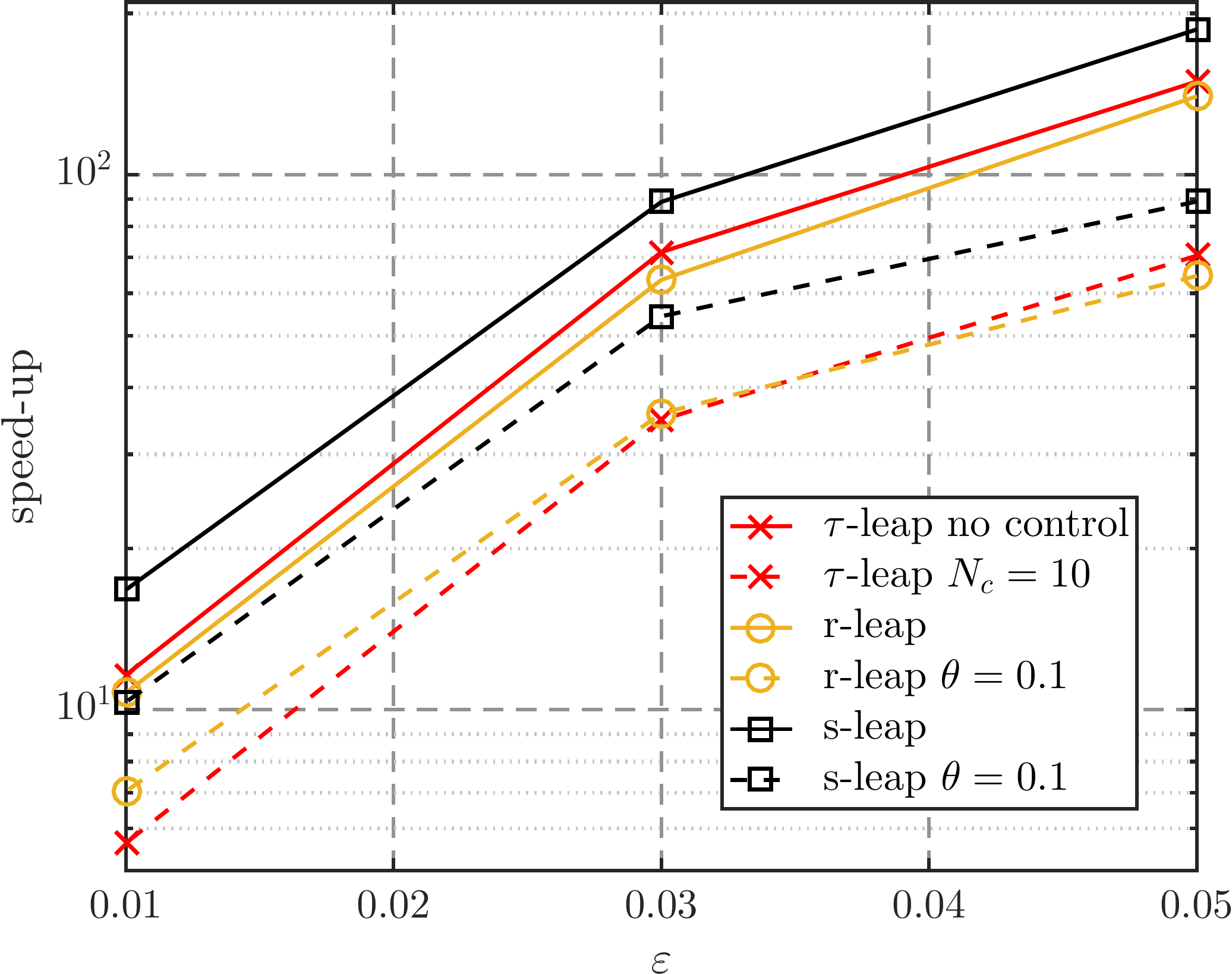}
\caption{Errors and efficiency for the LacZ/LacY system discussed in \cref{sec:lac}, with small initial population. }
\label{fig:lac:error:time}
\end{figure}

The system with small initial population is simulated until $T_{\textrm{end}}=2100$. Since none of the reversible reactions approached partial equilibrium during this time interval, only explicit methods are reported. To control the appearance of negative species, the $\tau$-leaping algorithm is used with control parameter $N_c=10$ \cite{Cao:2005}, while in the $S$-leaping and $R$-leaping we used $\theta=0.1$ as suggested in \cite{Auger:2006}. For comparison purposes, all three methods are also considered without the control mechanism. The frequency of reordering in the $R$-leaping and $S$-leaping is set to $p=10000$ as proposed in \cite{Auger:2006}. \cref{fig:lac:error:time} (right) shows the speed-up for the leap methods over SSA for $T_{\textrm{end}}=2100$. A single evaluation of the SSA for time $T_{\textrm{end}}=2100$ takes around $45$ min, making the evaluation of the models accuracy at this time point computationally expensive. Instead, \cref{fig:lac:error:time} (left) reports the error for all methods over the time interval $[0,100]$. For this system the error is averaged over the species TrLacZ2, TrRbsLacZ, and RbsribsomeLacY.

The $\tau$-leaping algorithm, as presented in \cref{alg:Tau}, executes mainly SSA steps and provides almost no speed up over SSA. Therefore, we turned off the SSA execution in the reported $\tau$-leaping algorithms. The leap methods without the control mechanism provide better speed-up over SSA, however their accuracy is reduced due to the high rejection rate. The sampling of reaction channels from the correlated binomial distribution in the $R$-leaping and $S$-leaping leads to lower rejection rate in comparison with the $\tau$-leaping, which is also reflected by the lower accuracy of the $\tau$-leaping method.

The control mechanisms in all leap methods results in high accuracy, at the cost of slightly reduced performance. The error reported in \cref{fig:lac:error:time} (left) is relatively constant and do not scale with $\varepsilon$, since the accuracy of these leap methods is mainly restricted by the mechanism preventing appearance of the negative species. The $S$-leaping reached comparable accuracy with the $R$-leaping, since they both use similar control mechanisms. On the other hand, the $\tau$-leaping considers most reactions critical and thus advance them with SSA, which lead to higher accuracy. The $R$-leaping and $S$-leaping algorithms benefit from the reordering of reaction channels and outperform the $\tau$-leaping. Moreover, since the stiffness of the system changes over time, the $S$-leaping outperforms both methods.

\begin{figure}
\centering
  	\includegraphics[width=0.45\textwidth]{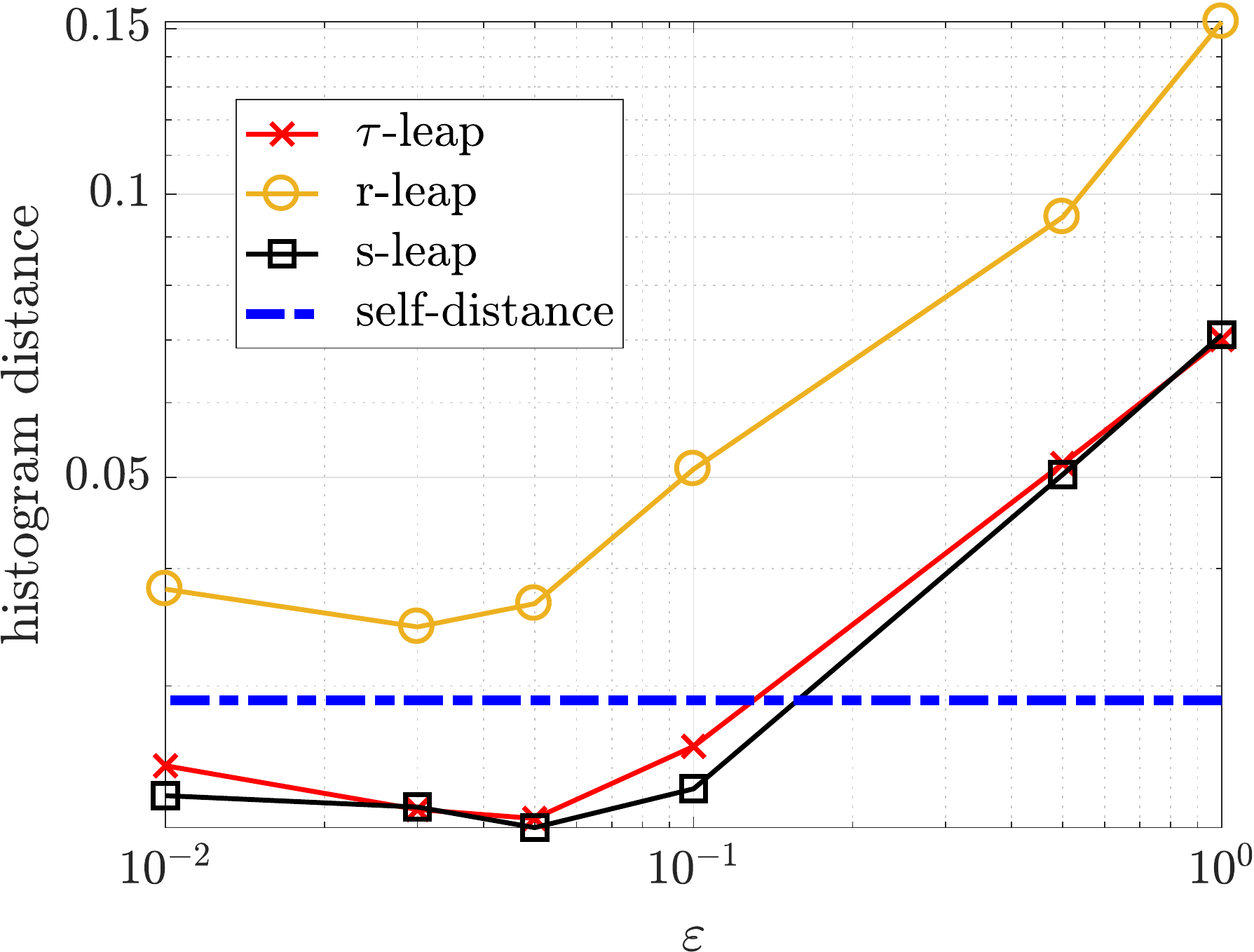}
	\includegraphics[width=0.45\textwidth]{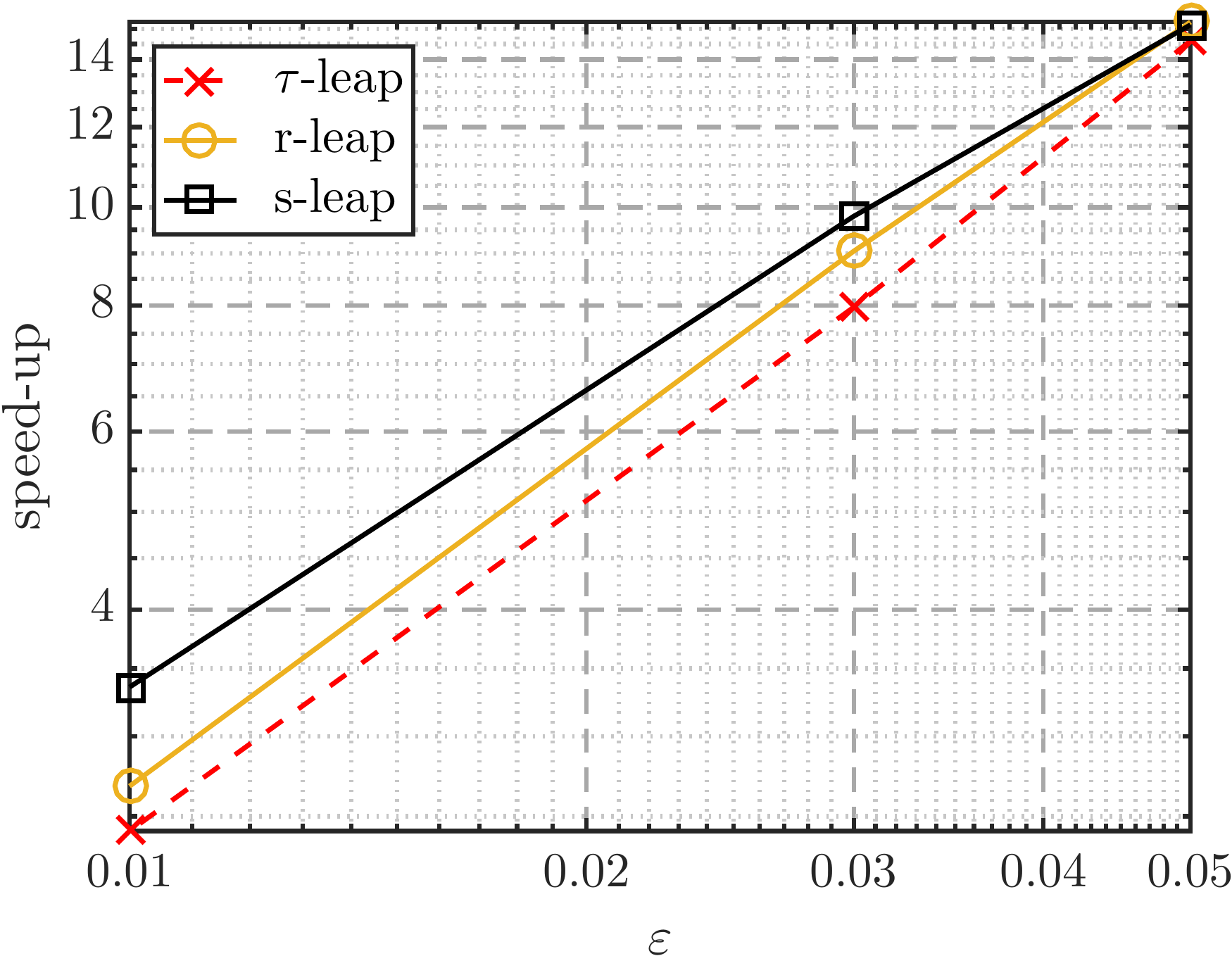}
\caption{Errors and efficiency for the LacZ/LacY system discussed in Section \ref{sec:lac}, with large initial population. }
\label{fig:lac:big:error:time}
\end{figure}

The system with big initial population is evolved until time $T_{\textrm{end}}=100$. As before, the $\tau$-leaping with the SSA steps performs mostly SSA and therefore the SSA step was disabled. Since all species appear in relatively large populations, the leap methods are considered without the control of negative population. The performance and accuracy of all methods is shown  in \cref{fig:lac:big:error:time}. The $S$-leaping algorithm again outperforms both the $\tau$-leaping and $R$-leaping method due to the combined advantages inherited from the both methods.
%
%
%
%~~~~~~~~~~~~~~~~~~~~~~~~~~~
\begin{table}
\begin{center}
\renewcommand*{\arraystretch}{1.2}
\begin{tabular}{ c | c | c }
& Reaction   & Reaction Rate \\ \hline
 $R_{1}$  &   PLac + RNAP $\longrightarrow$ PLacRNAP   &  0.17  \\  \hline
 $R_{2}$  &   PLacRNAP $\longrightarrow$ PLac + RNAP   &  10 \\  \hline
 $R_{3}$  &   PLacRNAP $\longrightarrow$ TrLacZ1   &  1 \\  \hline
 $R_{4}$  &   TrLacZ1 $\longrightarrow$ RbsLacZ + PLac + TrLacZ2   &  1 \\  \hline
 $R_{5}$  &   TrLacZ2 $\longrightarrow$ TrLacY1   &  0.015  \\  \hline
 $R_{6}$  &   TrLacY1 $\longrightarrow$ RbsLacY + TrLacY2   &  1  \\  \hline
 $R_{7}$  &   TrLacY2 $\longrightarrow$ RNAP   &  0.36  \\  \hline
 $R_{8}$  &   Ribosome + RbsLacZ $\longrightarrow$ RbsribosomeLacZ   &  0.17  \\  \hline
 $R_{9}$  &   RbsribosomeLacZ $\longrightarrow$ Ribosome + RbsLacZ   &  0.45  \\  \hline
 $R_{10}$  &   Ribosome + RbsLacY $\longrightarrow$ RbsribosomeLacY   &  0.17  \\  \hline
 $R_{11}$  &   RbsribosomeLacY $\longrightarrow$ Ribosome + RbsLacY   &  0.45  \\  \hline
 $R_{12}$  &   RbsribosomeLacZ $\longrightarrow$ TrRbsLacZ + RbsLacZ   &  0.4  \\  \hline
 $R_{13}$  &   RbsribosomeLacY $\longrightarrow$ TrRbsLacY + RbsLacY   &  0.4  \\  \hline
 $R_{14}$  &   TrRbsLacZ $\longrightarrow$ LacZ   &  0.015  \\  \hline
 $R_{15}$  &   TrRbsLacY $\longrightarrow$ LacY   &  0.036  \\  \hline
 $R_{16}$  &   LacZ $\longrightarrow$ dgrLacZ   &  6.42$\times 10^{-5}$  \\  \hline
 $R_{17}$  &   LacY $\longrightarrow$ dgrLacY   &  6.42$\times 10^{-5}$  \\  \hline
 $R_{18}$  &   RbsLacZ $\longrightarrow$ dgrRbsLacZ   &  0.3 \\  \hline
 $R_{19}$  &   RbsLacY $\longrightarrow$ dgrRbsLacY   &  0.3 \\  \hline
 $R_{20}$  &   LacZ + lactose $\longrightarrow$ LacZlactose   &  9.52 $\times10^{-5}$   \\  \hline
 $R_{21}$  &   LacZlactose $\longrightarrow$ product + LacZ   &  431  \\  \hline
 $R_{22}$  &   LacY $\longrightarrow$ lactose + LacY   &  14
\end{tabular}
\end{center}
\caption{The reaction network for the LacZ/LacY system discussed in Section \ref{sec:lac}.}
\label{table:lac}
\end{table}
%~~~~~~~~~~~~~~~~~~~~~~~~~~~
%
%
%
%
%~~~~~~~~~~~~~~~~~~~~~~~~~~~
%~~~~~~~~~~~~~~~~~~~~~~~~~~~
%
%
%~~~~~~~~~~~~~~~~~~~~~~~~~~~
%~~~~~~~~~~~~~~~~~~~~~~~~~~~

%==================================================

\newpage
%==================================================
%==========	Section: Conclusion ============
%==================================================
\section{Conclusion}
\label{sec:summary}
In this paper we have introduced the $S$-leaping,  an approximate algorithm for accelerating the SSA. The algorithm combines the advantages of two main approximate algorithms, the $\tau$-leaping and $R$-leaping.

The $S$-leaping method uses a time step selection, intrinsic to the $\tau$-leaping, which enables the extension of the algorithm to an implicit version. Furthermore, the $S$-leaping exploits the efficient sampling procedure from the $R$-leaping which reduce appearance of negative species. Moreover, the reordering of reaction channels inherited from the $R$-leaping, leads to a better performance of the $S$-leaping, compared to the $\tau$-leaping, in big and stiff systems. On the other hand, if a stiff system involves reversible reactions appearing close to equilibrium, then the implicit approach derived from the $\tau$-leaping accelerates the $S$-leaping by a few orders of magnitude in comparison to the explicit methods.

The performance of the proposed algorithm was tested on several examples, including a stiff, a non-stiff and a system involving slow and fast reactions with some species appearing in populations close to zero. In all test cases, accuracy of the $S$-leaping is similar to accuracy of the other accelerated methods. The performance of the $S$-leaping is comparable with the fastest method or even outperform both, the $\tau$-leaping and $R$-leaping methods. The $S$-leaping can be thus consider as optimal adaptive coupling of the $R$-leaping and $\tau$-leaping method.

Future work directions involve the extension of the $S$-leaping algorithm to systems with spatial component  by using compartment-based approach \cite{Erban:2009} or Brownian dynamics models \cite{Lipkova:2011} to extend the simulation framework for reaction-diffusion processes that arise in many biological systems.

\section{Acknowledgements}
PK and GA gratefully acknowledge support from the European Research Council (ERC) Advanced Investigator Award (No. 341117). {\color{black} The authors thank both anonymous reviewers for their insightful comments that have helped us to improve the content of the paper.}

%=================================================
\bibliographystyle{abbrv}
%\bibliography{SLeaping.bbl}
%\input{SLeaping.bbl}
%\bibliography{SLeaping_JCP}

\begin{thebibliography}{10}

\bibitem{Anderson:2011}
D.~F. Anderson and T.~G. Kurtz.
\newblock {\em Continuous Time Markov Chain Models for Chemical Reaction
  Networks}, pages 3--42.
\newblock Springer New York, New York, NY, 2011.

\bibitem{Auger:2006}
A.~Auger, P.~Chatelain, and P.~Koumoutsakos.
\newblock R-leaping: accelerating the stochastic simulation algorithm by
  reaction leaps.
\newblock {\em J Chem Phys}, 125(8):084103, Aug 2006.

\bibitem{Bayati:2011}
B.~Bayati, P.~Chatelain, and P.~Koumoutsakos.
\newblock {Adaptive mesh refinement for stochastic reaction-diffusion
  processes}.
\newblock {\em J. of Computational Physics}, 230(1):13--26, 2011.

\bibitem{Bayati:2010d}
B.~Bayati, H.~Owhadi, and P.~Koumoutsakos.
\newblock {A cutoff phenomenon in accelerated stochastic simulations of
  chemical kinetics via flow averaging (FLAVOR-SSA)}.
\newblock {\em Journal of Chemical Physics}, 133(24):1--7, 2010.

\bibitem{Cao1:2005}
Y.~Cao, D.~Gillespie, and L.~Petzold.
\newblock Multiscale stochastic simulation algorithm with stochastic partial
  equilibrium assumption for chemically reacting systems.
\newblock {\em Journal of Computational Physics}, 206(2):395--411, July 2005.

\bibitem{Cao:2006}
Y.~Cao, D.~Gillespie, and L.~Petzold.
\newblock Efficient step size selection for the tau-leaping simulation method.
\newblock {\em Journal of Chemical Physics}, 124(4):044109, Jan. 2006.

\bibitem{Cao:2005}
Y.~Cao, D.~T. Gillespie, and L.~R. Petzold.
\newblock Avoiding negative populations in explicit poisson tau-leaping.
\newblock {\em J Chem Phys}, 123(5):054104, Aug 2005.


\bibitem{Cao:2006b}
Y.~Cao, and L.~R. Petzold.
\newblock Accuracy limitations and the measurement of errors in the stochastic simulation of chemically reacting systems.
\newblock {\em J. of Computational Physics}, 212(1)6-24, 2006.


\bibitem{Cao:2007}
Y.~Cao, D.~T. Gillespie, and L.~R. Petzold.
\newblock Adaptive explicit-implicit tau-leaping method with automatic tau
  selection.
\newblock {\em J. of Chemical Physics}, 126(22):224101, June 2007.

\bibitem{Chattopadhyay:2013}
I.~Chattopadhyay, A.~Kuchina, G.~M. S{\"u}el, and H.~Lipson.
\newblock Inverse gillespie for inferring stochastic reaction mechanisms from
  intermittent samples.
\newblock {\em Proceedings of the National Academy of Sciences},
  110(32):12990--12995, 2013.

\bibitem{Erban:2007}
R.~Erban, J.~Chapman, and P.~Maini.
\newblock A practical guide to stochastic simulations of reaction-diffusion
  processes.
\newblock {\em https://arxiv.org/abs/0704.1908}, 2007.

\bibitem{Erban:2009}
R.~Erban and S.~J. Chapman.
\newblock Stochastic modelling of reaction--diffusion processes: algorithms for bimolecular reactions.
\newblock {\em Physical biology}, 6(4):046001, 2009.

\bibitem{Gibson:2000}
M.A.~Gibson and J.~Bruck
\newblock Efficient Exact Stochastic Simulation of Chemical Systems with Many Species and Many Channels
\newblock {\em The Journal of Physical Chemistry A}, 104(9):1876-1889, 2000.


\bibitem{Gillespie1976}
D.~T. Gillespie.
\newblock A general method for numerically simulating the stochastic time
  evolution of coupled chemical reactions.
\newblock {\em J. of Computational Physics}, 22(4):403--434, 1976.

\bibitem{Gillespie:1977}
D.~T. Gillespie.
\newblock Exact stochastic simulation of coupled chemical reactions.
\newblock {\em J. of physical chemistry}, 81(25):2340--2361, 1977.

\bibitem{Gillespie:2001}
D.~T. Gillespie.
\newblock Approximate accelerated stochastic simulation of chemically reacting
  systems.
\newblock {\em J. of Chemical Physics}, 115:1716, 2001.

\bibitem{Gillespie:2003}
D.~T. Gillespie and L.~R. Petzold.
\newblock Improved leap-size selection for accelerated stochastic simulation.
\newblock {\em J. of Chemical Physics}, 119:8229, 2003.

\bibitem{Kierzek:2002}
A.~M. Kierzek.
\newblock {STOCKS: STOChastic Kinetic Simulations of biochemical systems with
  Gillespie algorithm.}
\newblock {\em Bioinformatics (Oxford, England)}, 18(3):470--81, 2002.

\bibitem{Koumoutsakos:2013a}
P.~Koumoutsakos and J.~Feigelman.
\newblock {Multiscale stochastic simulations of chemical reactions with
  regulated scale separation}.
\newblock {\em J. of Computational Physics}, 244:290--297, 2013.

\bibitem{Lipkova:2011}
J.~Lipkova, K.~C. Zygalakis, S.~J. Chapman, and R.~Erban.
\newblock Analysis of brownian dynamics simulations of reversible bimolecular
  reactions.
\newblock {\em SIAM Journal On Applied Mathematics}, 71(3):714--730, 2011.

\bibitem{Maamar:2007}
H.~Maamar, A.~Raj, and D.~Dubnau.
\newblock Noise in gene expression determines cell fate in bacillus subtilis.
\newblock {\em Science}, 317(5837):526--529, 2007.

\bibitem{Mjolsness:2009}
E.~Mjolsness, D.~Orendorff, P.~Chatelain, and P.~Koumoutsakos.
\newblock An exact accelerated stochastic simulation algorithm.
\newblock {\em J. of Chemical Physics}, 130:144110, 2009.

\bibitem{Rathinam:2003}
M.~Rathinam, L.~R. Petzold, Y.~Cao, and D.~T. Gillespie.
\newblock Stiffness in stochastic chemically reacting systems: The implicit
  tau-leaping method.
\newblock {\em J. of Chemical Physics}, 119:12784, 2003.

\bibitem{Sandmann:2009}
W.~Sandmann.
\newblock Exposition and streamlined formulation of adaptive explicitimplicit
  tau-leaping.
\newblock Technical report, Citeseer, 2009.

\bibitem{Suel:2006}
G.~M. S{\"u}el, J.~Garcia-Ojalvo, L.~M. Liberman, and M.~B. Elowitz.
\newblock An excitable gene regulatory circuit induces transient cellular
  differentiation.
\newblock {\em Nature}, 440(7083):545--550, 2006.

\bibitem{Tian:2004}
T.~Tian and K.~Burrage.
\newblock Binomial leap methods for simulating stochastic chemical kinetics.
\newblock {\em The Journal of chemical physics}, 121:10356, 2004.



\end{thebibliography}

%==================================================

\end{document}